\documentclass{amsart}
\usepackage{amsaddr}
\usepackage[margin=2cm,letterpaper]{geometry}

\usepackage{amssymb,enumerate}
\usepackage{xr-hyper}
\newcommand{\csp}[1]{\begin{NoHyper}\cite[#1]{PavlovScholbach:Spectra}\end{NoHyper}}
\newcommand{\csy}[1]{\begin{NoHyper}\cite[#1]{PavlovScholbach:Symmetry}\end{NoHyper}}
\usepackage[hypertex,backref=section,linktoc=all]{hyperref}
\usepackage[all,cmtip]{xy}
\SilentMatrices
\usepackage[utf8]{inputenc} 

\def\typeout#1{}
\makeatletter
\@addtoreset{equation}{section}
\@addtoreset{equation}{subsection}
\makeatother

\theoremstyle{definition}
\newtheorem{Defi}[equation]{Definition} \newcommand{\defi}{\begin{Defi}} \newcommand{\xdefi}{\end{Defi}} \newcommand{\refde}[1]{Definition~\ref{defi--#1}}
\newtheorem{Bsp}[equation]{Example} \newcommand{\exam}{\begin{Bsp}} \newcommand{\xexam}{\end{Bsp}} 

\theoremstyle{remark}
\newtheorem{Bem}[equation]{Remark} \newcommand{\rema}{\begin{Bem}} \newcommand{\xrema}{\end{Bem}} \newcommand{\refre}[1]{Remark~\ref{rema--#1}}
\newtheorem{Nota}[equation]{Notation} \newcommand{\nota}{\begin{Nota}} \newcommand{\xnota}{\end{Nota}} 

\theoremstyle{plain}
\newtheorem{Theo}[equation]{Theorem} \newcommand{\theo}{\begin{Theo}} \newcommand{\xtheo}{\end{Theo}} \newcommand{\refth}[1]{Theorem~\ref{theo--#1}}
\newtheorem{Satz}[equation]{Proposition} \newcommand{\prop}{\begin{Satz}} \newcommand{\xprop}{\end{Satz}} \newcommand{\refpr}[1]{Proposition~\ref{prop--#1}}
\newtheorem{Lemm}[equation]{Lemma} \newcommand{\lemm}{\begin{Lemm}} \newcommand{\xlemm}{\end{Lemm}} \newcommand{\refle}[1]{Lemma~\ref{lemm--#1}}
\newtheorem{Coro}[equation]{Corollary} \newcommand{\coro}{\begin{Coro}} \newcommand{\xcoro}{\end{Coro}} \newcommand{\refco}[1]{Corollary~\ref{coro--#1}}

\newcommand{\refsect}[1]{\S\ref{sect--#1}}
\newcommand{\refchap}[1]{\S\ref{chap--#1}}
\newcommand{\refit}[1]{(\ref{item--#1})}
\newcommand{\refeq}[1]{(\ref{eqn--#1})}

\newcommand{\pf}{\begin{proof}} \newcommand{\xpf}{\end{proof}}

\makeatletter\@beginparpenalty=10000\makeatother
\let\ppar\endgraf 

\newcommand{\category}[1]{\mathbf{#1}}
\newcommand{\Set}{\category{Set}} 
\newcommand{\Cat}{\category{Cat}} 
\newcommand{\simpl}{\category s} 
\newcommand{\sSet}{\simpl\category{Set}} 
\newcommand{\ssSeq}{\category{(s)Seq}} 
\newcommand{\sSeq}{\category{sSeq}} 
\newcommand{\sColl}{\category{sColl}} 
\newcommand{\ssColl}{\category{(s)Coll}} 
\newcommand{\sOper}{\category{sOper}} 
\newcommand{\ssOper}{\category{(s)Oper}} 
\newcommand{\nsSeq}{\category{Seq}} 
\newcommand{\nsColl}{\category{Coll}} 
\newcommand{\nsOper}{\category{Oper}} 
\newcommand{\Pairs}{\category{Pairs}} 
\newcommand{\Top}{\category{Top}} 
\newcommand{\Sets}{\category{Sets}} 
\newcommand{\Mod}{\category{Mod}} 
\newcommand{\Alg}{\mathbf{Alg}} 
\newcommand{\ssTree}{\category{(s)Tree}} 
\newcommand{\Ch}{\category{Ch}} 
\newcommand{\Ar}{\category{Ar}} 
\newcommand{\Fun}{\category{Fun}} 

\newcommand{\proj}{\mathrm{pro}} 
\newcommand{\inje}{\mathrm{in}} 

\newcommand{\CO}{\mathrm{CO}} 
\newcommand{\cof}{\mathrm{C}} 
\def\dom{\mathop{\rm dom}} 
\def\codom{\mathop{\rm codom}} 
\newcommand{\fib}{\mathrm{F}} 
\newcommand{\we}{\mathrm{W}} 
\newcommand{\cell}{\mathrm{cell}} 
\def\mcr{{\rm Q}} 
\def\Obj{\mathop{\rm Obj}} 

\newcommand{\ws}[1]{\mathop{\rm cof}({#1})} 
\newcommand{\colim}{\operatornamewithlimits{colim}} 
\newcommand{\cofib}{\operatorname{cofib}} 
\def\id{{\rm id}} 
\def\op{{\rm op}} 
\def\Mor{\mathop{\rm Mor}\nolimits} 
\def\Hom{\mathop{\rm Hom}\nolimits} 
\def\IHom{\underline{\Hom}} 
\def\End{\mathop{\rm End}\nolimits} 
\def\Free{\mathrm{Free}} 
\def\Aut{\mathop{\rm Aut}\nolimits} 
\def\Map{\mathop{\rm Map}\nolimits} 
\def\val{\mathrm{val}} 
\def\Env{\mathop{\rm Env}\nolimits} 
\def\Ax{\Sigma} 

\mathchardef\pp"2403 
\mathchardef\ppdom"2400 

\def\bigpp{\mathop{\mathchoice{\dobigpp\Huge}{\dobigpp\Large}{\dobigpp\normalsize}{\dobigpp\small}}}
\mathchardef\bigppchar"1403
\def\dobigpp#1{\vcenter{#1\kern.2ex\hbox{$\bigppchar$}\kern.2ex}}
\mathchardef\bigppdomchar"1400
\def\dobigppdom#1{\vcenter{#1\kern.2ex\hbox{$\bigppdomchar$}\kern.2ex}}

\def\NN{{\bf N}} 
\def\Q{{\bf Q}} 

\def\N{{\rm N}} 

\newcommand{\Ob}{\mathrm{Ob}} 
\newcommand{\CatOp}{\mathrm{Cat}} 
\newcommand{\DiagOp}{\mathrm{Diag}} 
\newcommand{\OpOp}{\mathrm{Oper}} 
\newcommand{\ssOpOp}{\mathrm{(s)Oper}} 
\newcommand{\Comm}{\mathrm{Comm}} 
\newcommand{\Ass}{\mathrm{As}} 
\def\Ei{{\rm E_\infty}} 
\def\Ai{{\rm A_\infty}} 
\def\E{{\rm E}} 
\def\B{{\rm B}} 
\newcommand{\Y}{\mathcal{Y}} 
\newcommand{\C}{\mathcal{C}} 
\newcommand{\D}{\mathcal{D}} 
\def\V{{\cal V}} 
\newcommand{\cO}{\mathcal{O}}

\def\To#1#2{\mathbin{\mathop{\count0=#1 \loop\ifnum\count0>0 \smash-\mkern-7mu \advance\count0 -1 \repeat \mathord\rightarrow}\limits^{#2}}} 
\let\x\times

\let\t\otimes
\let\r\rightarrow

\def\matrix#1{\null\,\vcenter{\normalbaselines
    \ialign{\hfil$##$\hfil&&\quad\hfil$##$\hfil\crcr
      \mathstrut\crcr\noalign{\kern-\baselineskip}
      #1\crcr\mathstrut\crcr\noalign{\kern-\baselineskip}}}\,}
\def\vcd#1{\def\normalbaselines{\baselineskip20pt\lineskip1pt\lineskiplimit0pt } \harrowsize#1 \matrix}
\def\cd{\vcd{2em}}
\newdimen\harrowsize
\def\mapright#1{\smash{\mathop{\hbox to\harrowsize{\rightarrowfill}}\limits^{#1}}}

\def\mapdown#1{\Big\downarrow\rlap{$\vcenter{\hbox{$\scriptstyle#1$}}$}}

\catcode`@=11
\let\over\@@over
\let\atop\@@atop
\let\above\@@above
\let\overwithdelims\@@overwithdelims
\let\atopwithdelims\@@atopwithdelims
\let\abovewithdelims\@@abovewithdelims
\def\eqalign#1{\null\,\vcenter{\openup\jot\m@th
  \ialign{\strut\hfil$\displaystyle{##}$&$\displaystyle{{}##}$\hfil
      \crcr#1\crcr}}\,}
\newskip\xcentering
\def\eqalignno#1{\displ@y \tabskip\xcentering
  \halign to\displaywidth{\hfil$\@lign\displaystyle{##}$\tabskip\z@skip
    &$\@lign\displaystyle{{}##}$\hfil\tabskip\xcentering
    &\llap{$\@lign##$}\tabskip\z@skip\crcr
    #1\crcr}}
\def\cases#1{\left\{\,\vcenter{\normalbaselines\m@th
    \ialign{$##\hfil$&\quad##\hfil\crcr#1\crcr}}\right.}
\def\@writetocindents{}
\def\eqlabel#1{\refstepcounter{equation}\label{eqn--#1}\ifmmode\ifinner\else\eqno\fi\fi\hbox{\@eqnnum}} 
\def\cal{\fam\tw@}

\catcode`@=12

\begin{document}

\title{Admissibility and rectification of colored symmetric operads}

\def\obfuscate#1#2{\rlap{\hphantom{#2}@#1}#2\hphantom{@#1}}

\author{Dmitri Pavlov}
\address{Faculty of Mathematics, University of Regensburg\\
Department of Mathematics and Statistics, Texas Tech University}
\email{https://dmitripavlov.org/}

\author{Jakob Scholbach}
\address{Mathematical Institute, University of M\"unster}
\email{https://wwwmath.uni-muenster.de/u/jakob.scholbach/}

\subjclass[2010]{Primary 55P48, 18D50; Secondary 55U35, 55P43, 18G55, 18D20.}


\centerline{\font\tfont=cmss17 \tfont\shorttitle}

\bigskip
{\tabskip0pt plus 1fil
\halign to\hsize{&\hfil#\hfil\cr
{\bf Dmitri Pavlov\/}\cr
Faculty of Mathematics, University of Regensburg\cr
Department of Mathematics and Statistics, Texas Tech University\cr
\href{https://dmitripavlov.org/}{https:/\negthinspace/dmitripavlov.org/}\vadjust{\medbreak}\cr
{\bf Jakob Scholbach\/}\cr
Mathematical Institute, University of M\"unster\cr
\href{https://wwwmath.uni-muenster.de/u/jakob.scholbach/}{https:/\negthinspace/wwwmath.uni-muenster.de/u/jakob.scholbach/}\cr
}}

\begin{abstract}
We establish a highly flexible condition that guarantees that all colored symmetric operads in a symmetric monoidal model category are admissible,
i.e., the category of algebras over any operad admits a model structure transferred from the original model category.
We also give a necessary and sufficient criterion that ensures that a given weak equivalence of admissible operads admits rectification,
i.e., the corresponding Quillen adjunction between the categories of algebras is a Quillen equivalence.
In addition, we show that Quillen equivalences of underlying symmetric monoidal model categories yield
Quillen equivalences of model categories of algebras over operads.
Applications of these results include enriched categories, colored operads, prefactorization algebras, and commutative symmetric ring spectra.
\end{abstract}

\makeatletter\@setabstract\makeatother

\setcounter{tocdepth}{1}
\tableofcontents

\numberwithin{equation}{section}

\section{Introduction}

This paper is devoted to the model-categorical study of operads and their algebras.
The concept of an algebra over a colored symmetric operad allows for a uniform treatment of algebraic structures that produce an output from multiple inputs,
subject to some symmetry constraints.
For example, a commutative monoid~$X$ in a symmetric monoidal category~$\C$ is specified by $\Sigma_n$-equivariant maps $X^{\t n} \r X$,
subject to the usual associativity and unitality constraints.
In a seemingly artificial way, this can be rewritten as
$$\Comm_n \t_{\Sigma_n} X^{\t n} \r X,$$
where $\Comm$ is the so-called commutative operad, which satisfies $\Comm_n = 1$, the monoidal unit.
More generally, an algebra of a single-colored operad $O$ is an object $A \in \C$ together with maps
$$O_n \t_{\Sigma_n} A^{\t n} \r A,$$
which are compatible with the multiplication in $O$ in a suitable sense.
Colored symmetric operads, also known as symmetric multicategories, are a many-objects version of ordinary operads.
They allow input from more than one object.
For example, there is a two-colored operad whose algebras are pairs $(R, M)$, where $R$ is a commutative monoid in $\C$ and $M$ is an $R$-module.
Interestingly, operads themselves are algebras over a certain operad.

Symmetric operads and their algebras, which were first introduced by May, are ubiquitous in homotopy theory and beyond.
A prototypical example is the $m$-fold loop space $\Omega^mX$ of some topological space~$X$: concatenation of paths yields a multiplication map
$$\mu_n\colon(\Omega^mX)^n \r \Omega^mX,$$
which is neither associative nor commutative, but only associative and commutative up to homotopy.
This and the compatibility of these homotopies for various $n$ is concisely encoded in the fact that $\Omega^m X$ is an algebra over some operad~$O$,
meaning that there are maps (for all~$n$, and compatible with each other):
$$O_n \x_{\Sigma_n} (\Omega^mX)^n \r \Omega^mX.$$
If $O_n$ was just a point, then this would mean that the multiplication on $\Omega^mX$ is strictly commutative and associative, which it is not.
However, $O$ can be chosen to be the little disks operad~$\E_m$.
For $m=\infty$ these levels~$O_n$ are contractible spaces,
which can be interpreted as saying that infinite loop spaces are homotopy coherent commutative monoids.
Recently, $\E_n$-algebras have been attracting a lot of attention in questions related to factorization homology (also known as topological chiral homology)
and Goodwillie calculus of functors.

Our first main theorem is a highly flexible existence criterion for a model structure on algebras over operads in a model category.
This is a powerful tool for homotopical computations related to algebras over operads, such as the loop space.

\theo \label{theo--admissible} (See Theorems \ref{theo--O.Alg},
\ref{prop--cofibrant.strongly.admissible},
\ref{theo--O.strongly.admissible}.)
Suppose that $\C$ is a symmetric monoidal model category that is \emph{symmetric h-monoidal} and satisfies some minor technical assumptions.
Then any symmetric $W$-colored operad $O$ is \emph{admissible},
i.e., the category $\Alg_O(\C)$ of $O$-algebras carries a model structure whose weak equivalences and fibrations are inherited from~$\C$.

The forgetful functor $\Alg_O(\C) \r \C^W$ preserves cofibrant objects and cofibrations between them if $\C$ is \emph{symmetroidal} and $O$ is weakly well-pointed (essentially, this means $O$ is levelwise injectively cofibrant).
Alternatively, the same statement is true if $O$ is well-pointed (essentially, $O$ is levelwise projectively cofibrant).
\xtheo

This admissibility result is widely applicable because its assumptions are satisfied for many basic model categories such as simplicial sets, topological spaces, simplicial presheaves,
chain complexes of rational vector spaces.
It does not apply to chain complexes of abelian groups, and in fact the commutative operad is provably not admissible in this category.
Moreover, as was shown in~\cite{PavlovScholbach:Symmetry}, symmetric h-monoidality (and similarly with symmetroidality and symmetric flatness)
are stable under transfer and monoidal left Bousfield localizations, which allows to easily promote these properties from basic model categories to more advanced model categories,
such as spectra.
The latter are shown in \cite{PavlovScholbach:Spectra} to be symmetric h-monoidal, symmetroidal, and symmetric flat.

The key condition of symmetric h-monoidality is a symmetric strengthening of the h-monoidality condition.
The latter was introduced by Batanin and Berger in~\cite{BataninBerger:Homotopy} and is closely related to the monoid axiom.
Essentially, it means that for any object~$Y$ in $\Sigma_n \C$ (objects of~$\C$ with a $\Sigma_n$-action) and any cofibration~$f$, the map
$$Y \t_{\Sigma_n} s^{\pp n} := (Y \t s^{\pp n})_{\Sigma_n}$$
is an h-cofibration, which is a weak equivalence if $f$ is an acyclic cofibration.
Here $f^{\pp n}$ is the $n$-fold pushout product of~$f$.
Symmetroidality is a related condition, obtained by replacing ``h-cofibration'' above by ``cofibration'' and $Y \t -$ by $y \pp -$ for some map~$y$.

In practice, a frequent question is how to replace algebras over some operad by those over a weakly equivalent operad.
A common class of (weakly equivalent) operads is defined by the condition that $O_n$ is a contractible space and has a free $\Sigma_n$-action.
Such operads are referred to as $\Ei$-operads, and many variants of little disks operads are of this type.
One can therefore ask whether $\Omega^\infty X$, together with the multiplications $\mu_n$, is weakly equivalent to some space with a strictly commutative and associative multiplication.
In this example, it is well known that connected $\Ei$-spaces with nontrivial Postnikov invariants, e.g., the identity component of the space~$\Omega^\infty\Sigma^\infty S^0$,
\emph{cannot} be strictified to a simplicial abelian group.
Indeed by a classical result of Moore~\cite[Theorem~3.29]{Moore:Semi}, connected simplicial abelian groups have trivial Postnikov invariants.

The following rectification theorem identifies a criterion when a rectification of operadic algebras is possible.

\theo \label{theo--Quillen.invariance} (See \refth{rect.colored.operad}.)
For any map of admissible operads $O \r P$ in a symmetric monoidal model category, there is a Quillen adjunction
$$\Alg_O (\C) \rightleftarrows \Alg_P(\C).$$
Provided that $\C$ satisfies some minor technical assumptions, it is a Quillen equivalence if and only if $O \r P$ is a \emph{symmetric flat} map in $\C$.
\xtheo

The symmetric flatness condition essentially requires that the map
$$O_n \t_{\Sigma_n} X^{\t n} \r P_n \t_{\Sigma_n} X^{\t n}$$
is a weak equivalence for all cofibrant objects $X$ and all $n \ge 0$.
If $\C$ is the model category of rational chain complexes, this condition holds for all weak equivalences $O \r P$.
In \cite{PavlovScholbach:Spectra}, we show that the same is true for symmetric spectra in an abstract model category.
However, this condition does not hold for all maps in simplicial sets, in particular, it fails for the components of $\Ei\to\Comm$.
This matches the above observation of the nonrectifiability of $\Ei$-algebras to strictly commutative simplicial monoids.
Nevertheless, it is satisfied for any pair of $\Ei$~operads in simplicial sets, which shows that the algebras over such operads are all Quillen equivalent to each other.

As a consequence of this rectification result, we obtain \refth{quasicategorical.rectification}, which relates algebras over operads in the strict sense, as above,
and algebras over quasicategorical operads as introduced by Lurie.

Operads and their algebras in different model categories also behave as nicely as possible.
Such a result allows to replace $\C$ by a more convenient model category, which is often necessary in practice.

\theo (See \refth{transport}.)
For any Quillen equivalence
$$F \colon \C \rightleftarrows \D\colon G$$
between symmetric monoidal model categories as above, where $F$ is symmetric oplax monoidal such that
the canonical maps $F\mcr (1_\C) \to 1_\D$ and $F(C \t C') \to F(C) \t F(C')$ are weak equivalences for all cofibrant objects $C, C' \in \C$,
there is a Quillen equivalence of the categories of $W$-colored (symmetric) operads
$$F^\ssOper : \ssOper(\C) \rightleftarrows \ssOper(\C') : G.$$
Moreover, there is a Quillen equivalence for any cofibrant (symmetric) operad~$O$,
$$F^\Alg: \Alg_{O} (\C) \rightleftarrows \Alg_{F^\ssOper(O)} (\D): G.$$
\xtheo

The admissibility and rectification of nonsymmetric and symmetric operads is a topic that was addressed by various authors.
Spitzweck has shown the existence of a semi-model structure for special symmetric operads,
namely those whose underlying symmetric sequence is projectively cofibrant (which roughly means that $\Sigma_n$ acts freely on $O_n$) \cite[Theorem~4.7]{Spitzweck:Operads}.
This rules out the commutative operad, whose algebras are commutative monoids.
The admissibility of the commutative operad was shown by Lurie under the assumption of symmetroidality of the commutative operad
(see \cite[Lemma~4.5.4.11.(1) and Proposition~4.5.4.6]{Lurie:HA}).
An independent account of this result was later given by White~\cite[Theorem~3.2]{White:Model}.
The admissibility of all operads was shown by Elmendorf and Mandell for $\C = \sSet$ \cite[Theorem~1.3]{ElmendorfMandell:Rings},
Berger and Moerdijk \cite{BergerMoerdijk:Axiomatic}, and Caviglia \cite{Caviglia:Model} for colored operads.
The latter two results use an assumption on the path object, which serves to cut short a certain homotopical analysis of pushouts performed in this paper.
The path object argument was also used by Johnson and Yau to establish a model structure on colored PROPs \cite{JohnsonYau:Homotopy}.
PROPs are more general than symmetric operads in that not only multiple inputs, but also multiple outputs are allowed.
Harper showed the admissibility of all symmetric operads in simplicial symmetric spectra~\cite{Harper:Symmetric}.
This was generalized by Hornbostel to spectra in simplicial presheaves~\cite{Hornbostel:Preorientations}.
Finally, Muro has shown the admissibility of all nonsymmetric operads~\cite{Muro:Homotopy, Muro:Corrections}.
A more detailed review of these results is found in \refsect{operad.admissible}.

Harper also established a rectification result under the assumption that every symmetric sequence is projectively cofibrant~\cite[Theorem~1.4]{Harper:Monoidal}.
This strong assumption applies to categories such as rational chain complexes.
In this case, rectification is due to Hinich~\cite{Hinich:Homological}.
Lurie \cite{Lurie:HA} established rectification of $\Ei$-algebras in the context of $\infty$-operads,
again under a strong assumption that only applies to special model categories such as rational chain complexes.
These and further results are reviewed in \refsect{rectification}.

Thus, all previous results have either restrictions on the operad and/or on the category in which the operad lives.
Our results are applicable to all operads and to a very broad range of model categories.
This wide applicability results from the fact that conditions of symmetric h-monoidality, symmetroidality, and symmetric flatness
occurring above are stable under transfer and monoidal left Bousfield localization.
Thus, they are easily promoted from simplicial sets to simplicial presheaves, say.

In \refsect{symmetricity}, we recall the \emph{symmetricity properties} introduced in \cite{PavlovScholbach:Symmetry}: symmetric h-monoidality, symmetroidality, and symmetric flatness, and a few other basic notions on model categories.
As was shown in \csy{\ref{prop--transfer.monoidal},
\ref{theo--transfer.symmetric.monoidal}, \ref{prop--monoidal.localization}, \ref{theo--symmetric.monoidal.localization}}, these properties are stable transfer and monoidal Bousfield localizations.
Given that these two methods are the most commonly used tools to construct model structures, the admissibility and rectification results in this paper are applicable to a wide range of model categories.

In \refsect{colored.collections}, we start with a brief review of colored symmetric collections and the substitution product.
Symmetric operads are defined as monoids in this category $\sColl_W(\C)$.
In \refsect{operad.admissible}, we show that symmetric h-monoidality is the key condition needed to ensure the admissibility of arbitrary symmetric operads $O$, i.e., the existence of the transferred model structure on $O$-algebras.
In \refsect{strong.admissibility}, we show that symmetroidality is needed to additionally guarantee the strong admissibility of $O$, i.e.,
the functor forgetting the $O$-algebra structure preserves cofibrations with cofibrant source.
In \refsect{rectification}, we show the rectification of algebras of weakly equivalent symmetric operads.
In \refsect{transport}, we establish Quillen equivalences of operads and their algebras in different model categories.

We obtain the above-mentioned theorems by systematically using the symmetricity properties above
combined with Berger--Moerdijk and Spitzweck's description of certain pushouts of operads \cite{Spitzweck:Operads, BergerMoerdijk:Derived}.
In \refchap{applications}, we finish this paper with examples and applications ranging from low-dimensional category theory to prefactorization algebras.

We thank the referee for suggesting \refle{projectively.cofibrant.symmetric.flat}.
We thank Clemens Berger, Giovanni Caviglia, Denis-Charles Cisinski, John Harper, Jacob Lurie, Birgit Richter, Brooke Shipley, and David White for helpful conversations.
We thank Thomas Nikolaus for a discussion that led to \refth{quasicategorical.rectification}.
We thank Yonatan Harpaz, Joost Nuiten, and Matan Prasma for offering a simple fix for a gap in the proof of \refpr{sifted.homotopy.colimits}.
We thank Miguel Barrero for pointing out a gap in the previous proof of \refth{rect.colored.operad}.
This work was partially supported by the SFB 878 grant.

\section{Symmetricity properties}
\label{sect--symmetricity}

Let $\C$ be a symmetric monoidal model category in the sense of \cite[Definitions~4.1.6, 4.2.6]{Hovey:Model}, except that we do not require the unit axiom.
In this section, we briefly recall the symmetricity properties from \csy{\refchap{symmetric}},
which are the key conditions in the admissibility, strong admissibility, and rectification results of this paper (see Theorems \ref{theo--O.Alg}, \ref{theo--O.strongly.admissible}, \ref{theo--rect.colored.operad}).

We use the notation of \csy{especially \refsect{arrows}, \refde{multi}}.
In particular, in the definitions below, $n = (n_1, \ldots, n_e)$ is an arbitrary finite multi-index.
For a family $s = (s_1, \ldots, s_e)$ of maps in $\C$, $\Sigma_n := \prod_i \Sigma_{n_i}$ acts on
the pushout product $s^{\pp n} := \bigpp_i s_i^{\pp n_i}$.
A subscript $\Sigma_n$ denotes the coinvariants of the $\Sigma_n$-action, such as $- \t_{\Sigma_n} -$.

The concept of h-monoidality in Part \refit{h.monoidal} is due to Batanin and Berger \cite[Definition~1.7]{BataninBerger:Homotopy}.
Recall from the same work that an \emph{h-cofibration}
$f\colon X \r Y$ is a map such that in any pushout diagram
$$\xymatrix{
X \ar[d]_f \ar[r] & A \ar[d] \ar[r]^g & B \ar[d]\cr
X' \ar[r] & A' \ar[r]^{g'} & B',\cr
}$$
the map~$g'$ is a weak equivalence if $g$ is one.
If, in addition, $f$ is a weak equivalence, it is an \emph{acyclic h-cofibration}.

\defi
\label{defi--summary.symmetricity}
Suppose $\C$ is a symmetric monoidal model category.
\begin{enumerate}[(i)]
\item
$\C$ is \emph{admissibly generated} if it is cofibrantly generated and if the (co)domains of a set $I$ of generating cofibrations (equivalently, by \cite[Corollary~10.4.9]{Hirschhorn:Model}, all cofibrant objects) are small with respect to the subcategory
$$\cell (Y\otimes_{\Sigma_n} s^{\pp n})\eqlabel{cellular.algebraic}$$
for any finite family~$s$ of cofibrations, and any object~$Y \in \Sigma_n \C$.
As usual, $\cell$ denotes the closure of a class of maps under pushouts and transfinite composition.
\item
\label{item--strongly.admissibly.generated}
$\C$ is \emph{strongly admissibly generated} if it is cofibrantly generated and if $\text{(co)dom}(I)$ are ($\aleph_0$-)compact (also known as finite) relative to \refeq{cellular.algebraic} \cite[Definition~10.8.1]{Hirschhorn:Model}.
\item
\label{item--h.monoidal}
$\C$ is \emph{h-monoidal} if the map $Y \t s$ is an (acyclic) h-cofibration for any (acyclic) cofibration $s$, and any object~$Y \in \C$.
\item
$\C$ is \emph{symmetric h-monoidal} if $Y\otimes_{\Sigma_n} s^{\pp n}$ is an (acyclic) h-cofibration for any finite family~$s$ of (acyclic) cofibrations, and any~$Y \in \Sigma_n \C$.
\item
Let $\Y=(\Y_n)_{n\ge1}$ be a collection of classes $\Y_n$ of morphisms in $\Sigma_n\C$, where $n \ge 1$ is any finite multi-index.
We suppose that for $y \in \Y_n$, $y \pp -$ preserves \emph{injective} (acyclic) cofibrations in $\Sigma_n\C$, i.e., those maps that are (acyclic) cofibrations in $\C$.
Then $\C$ is \emph{$\Y$-symmetroidal} if the morphism
$$y \pp_{\Sigma_n} s^{\pp n}$$
is an (acyclic) cofibration in~$\C$ for all finite families $s$ of (acyclic) cofibrations and all maps $y \in \Y_n$.
If $\Y_n = \cof_{\Sigma_n^\inje \C}$ (injective cofibrations), we say that $\C$ is \emph{(acyclic) symmetroidal}.
\item
A weak equivalence $y$ is flat if $y \pp s$ is a weak equivalence in $\C$ for any cofibration $s$.
$\C$~is \emph{flat} if all weak equivalences are flat.
\item
A weak equivalence $y \in \Sigma_n \C$ is called \emph{symmetric flat} if $y \pp_{\Sigma_n} s^{\pp n}$ is a weak equivalence (in $\C$) for any family $s$ of cofibrations.
We say that $\C$ is symmetric flat if all weak equivalences in $\Sigma_n \C$ are so.
\end{enumerate}
\xdefi

These conditions are usually stable under weak saturation (i.e., saturation under pushouts, transfinite compositions, and retracts).
Thus, they only have to be checked for \emph{generating} (acyclic) cofibrations $s$.
Simplicial sets with their standard model structure are symmetroidal, symmetric h-monoidal, and flat (but not symmetric flat).
The same is true for simplicial presheaves with the projective, injective, or local (with respect to some topology) model structures, and also for simplicial modules.

For any commutative ring $R$, chain complexes of $R$-modules with their projective model structure are flat and h-monoidal.
They are symmetroidal, symmetric h-monoidal, and symmetric flat if and only if $R$ contains~$\Q$.

The admissible generation is automatic if $\C$ is combinatorial \cite[Definition~A.2.6.1]{Lurie:HTT}.
Moreover, topological spaces are admissibly generated, symmetric h-monoidal, and symmetroidal.

To check symmetricity properties of more involved model categories, one can use the fact that the properties above are stable under transfer (appropriately compatible with the monoidal structure), and monoidal Bousfield localizations.
Combining these principles, we show in \csp{\refsect{symmetricity.2}} that spectra with values in a flat, h-monoidal (but not necessarily symmetric flat nor symmetric h-monoidal) category $\C$, with the positive stable model structure, are symmetric flat, symmetroidal, and symmetric h-monoidal.
In particular, this allows to replace $\C$ by a Quillen equivalent symmetric flat and symmetric h-monoidal model category.

The reader is referred to \csy{Theorems \ref{theo--symmetric.weakly.saturated}, \ref{theo--transfer.symmetric.monoidal}, \ref{theo--symmetric.monoidal.localization}, \refchap{examples}} for precise statements of the above facts and further examples.

Many results below include a condition that weak equivalences in $\C$ are stable under transfinite compositions or filtered colimits.
This condition is satisfied if $\C$ is cofibrantly generated and its generating cofibrations $I$ have compact domain and codomain or, slightly more generally, if $\C$ is pretty small in the sense of \csy{\refde{pretty.small}}.
This condition is satisfied for $\sSet$, $\Ch(\Mod_R)$, and many other basic model categories, but not for $\Top$.
However, $\Top$ is strongly admissibly generated, which is enough to conclude that the filtered colimits of the weak equivalences that actually occur (as a result of a cellular presentation of cofibrant objects) are indeed again weak equivalences.
We call $\C$ {\it tractable\/} if its (acyclic) cofibrations are contained in the weak saturation of (acyclic) cofibrations with cofibrant source (and target). This condition was introduced by Barwick \cite[Definition~1.21]{Barwick:Left}, who also includes the condition that  $\C$ be combinatorial.
Again, this holds for $\sSet$, $\Top,$ $\Ch(\Mod_R)$.
All three conditions are stable under monoidal Bousfield localizations and transfer, turning them into viable and effectively checkable conditions.

\section{Colored collections}
\label{sect--colored.collections}

In \S\ref{sect--colored.collections}--\ref{sect--enveloping.operad}, let $\C$ be a closed symmetric monoidal category.
In this section, we give a very brief overview of $W$-colored (symmetric) operads and colored modules over them (e.g., algebras over operads).
The reader can consult Gambino and Joyal~\cite{GambinoJoyal:Operads} for more details.
Constructions in this section involve a set~$W$, whose elements are called {\it colors}.
The reader may assume that $W$ has exactly one element, which yields ordinary operads.

$W$-colored symmetric operads in~$\C$ are defined as monoids in a certain monoidal category~$(\sColl_W(\C),\circ)$
and $V$-colored modules over a given $W$-colored (symmetric) operad~$O$ are defined as left modules over~$O$ in the category~$(\sColl_{V,W}(\C),\circ)$,
which itself is a left module over the monoidal category~$(\sColl_W(\C),\circ)$.
The idea behind~$(\sColl_W(\C),\circ)$ is that an object in~$\sColl_W(\C)$ encodes all possible operations,
whereas the monoid structure encodes the composition of operations.
Operations have a multisource, consisting of a finite family of colors,
and a target, which is a single color.
Furthermore, for any operation we can permute elements in its source and obtain another operation.
Operations with a fixed multisource and target form an object of~$\sColl_W(\C)$.
Likewise, an object in~$\sColl_{V,W}(\C)$ encodes operands that can be acted upon from the left by operations in a $W$-colored operad
and the left module structure encodes these actions.
The operands are encoded by a $V$-valued multisource and a target in~$W$.
Thus, the data of all operations can be encoded as a $\C$-valued presheaf on a certain groupoid~$\sSeq_W$ or~$\sSeq_{V,W}$, which we define first.

We simultaneously treat symmetric and nonsymmetric $W$-colored operads
with values in a symmetric monoidal category~$\C$, indicating the modifications necessary for the symmetric case in parentheses.
That is, we write $\ssOper$ to mean either $\sOper$ (symmetric operads) or $\nsOper$ (nonsymmetric operads) and likewise for other notations.

\defi
Given two sets~$V$, $W$, define the groupoid of \emph{(symmetric) $V,W$-sequences} as
$$\ssSeq_{V,W} := \ssSeq_V^\x\times W,$$
where $W$ denotes a category with objects~$W$ and identities as morphisms
and $\ssSeq_V^\x$ is the category of functions $s\colon I\to V$, where $I$ is a finite ordered set (respectively, finite unordered set, in the symmetric case) set
and morphisms $s\to s'$ are isomorphisms of ordered (respectively, unordered) sets
$f\colon I\to I'$ such that~$s=s'f$.
We abbreviate $\ssSeq_W :=\ssSeq_{W,W}$.
\xdefi

The idea is that an object~$(s,t)$ in~$\ssSeq_W^\x\times W$ encodes multisource~$s$ and target~$t \in W$.
Morphisms in~$\sSeq_W^\x$ account for the fact that one can permute sources in the symmetric case.
In the nonsymmetric variant $\nsSeq_W^s$, no permutation of multisources is allowed.
If $W=\{*\}$, then $\ssSeq_W$ is the category $\NN$ of finite ordered sets and identity morphisms (respectively, the category $\Sigma$ of symmetric sequences, i.e., finite sets and bijections).
Their objects can be interpreted as arities.
For some $s:I \r W$, we write $\Ax_s := \Aut_{\ssSeq_W^\x}(s)$.
In the nonsymmetric case, this group is trivial.
In the symmetric case, there is an isomorphism
$$\Ax_s = \prod_{w \in W} \Sigma_{s^{-1}(w)}.
\eqlabel{Ax.s}$$
For example, if $W = \{*\}$, then $\Ax_s = \Sigma_{\sharp I}$.

Given a (symmetric) sequence $X \in \ssSeq_W$, we write $X_0 \in \C^W$ for the restriction to objects with empty multisource, i.e., $s\colon \emptyset \r W$.
We refer to this by saying that $X_0$ is concentrated in \emph{degree~$0$}.
We refer to the $X_{s,w}$ with $s\colon I \r W$ satisfying $\sharp I = 1$, $s(i) = w$ as the \emph{unit degrees} and will write $X_{w,w}$ in this case.
The remaining components are called the \emph{nonunit degrees}.

\defi
Given symmetric monoidal categories $\V$~and~$\C$ such that $\C$ is enriched over~$\V$,
for a given pair of sets $V$~and~$W$ define the categories
$$\ssColl_{V,W}(\C) := \Fun(\ssSeq_{V,W}^\op,\C),$$
where $\Fun$ denotes the $\V$-enriched category of functors.
Set $$\ssColl_W(\C)=\ssColl_{W,W}(\C),$$ which we call the category of {\it $W$-colored (symmetric) collections in~$\C$}.
The category $\ssColl_W(\C)$ is a monoidal category and the category $\ssColl_{V,W}(\C)$ is a left module over~$\ssColl_W(\C)$
via the {\it substitution product\/}
$$\circ\colon\ssColl_{V,W}(\C)\times\ssColl_{U,V}(\C)\to\ssColl_{U,W}(\C).\eqlabel{substitution.product}$$
The substitution product of $F\in\sColl_{V,W}(\C)$ and $G\in\sColl_{U,V}(\C)$ can be computed as
the left Kan extension
$$\xymatrix{
T_{U,V,W} \ar[rr]^{F*G} \ar[d]_{\text{proj.}} & & \C \cr
\ssSeq_U \x W \ar[urr]_{F \circ G},}$$
where $T_{U,V,W}$ is the category whose objects are quadruples $(u\colon I\to U,v\colon J\to V,w\colon1\to W,f\colon I\to J)$, where $I$~and~$J$ are finite sets,
and morphisms are commutative diagrams
$$\xymatrix@=1.2pc{
& I \ar[dd]_i^\cong \ar[dl]_u \ar[rr]^f & & J \ar[dd]^j_\cong \ar[dr]^v \cr
U & & & & V \cr
& I' \ar[ul]^{u'} \ar[rr]_{f'} & & J' \ar[ur]_{v'},\cr
}$$
where $i$ and $j$ are isomorphisms and $w=w'$.
The functor $F*G$ sends an object~$(u,v,w,f)$ to $F(v,w)\otimes\bigotimes_{p\in J}G(u|_{f^{-1}(p)},p)$
and a morphism~$(i,j)$ to the isomorphism $F(j)\otimes\bigotimes_{p\in J}G(i|_{f^{-1}(p)})$.

The monoidal unit of $\ssColl_W$ is the $W$-colored collection that assigns the monoidal unit~$1\in\C$ to all unit degrees $(w,w)$, $w \in W$ and the initial object of~$\C$ to anything else.
We denote it by $1[1]$.
\xdefi

See Gambino and Joyal \cite[Theorem~10.2 and Remark~11.7]{GambinoJoyal:Operads} for additional details.
In the notation of Gambino and Joyal, $\cal R$ stands for~$\C$.

\exam
For example, for $U=\emptyset$, which is the special case relevant for algebras over colored operads,
$$(F*G)(v,w)=F(v,w)\otimes\bigotimes_{p\in J}G(p)$$
and $(F*G)(j)=F(j)\otimes\id$.

In the case $W=\{*\}$, the substitution product in $\sColl$ can be expressed concisely using the symmetric smash product~$\otimes$ on symmetric sequences
(see Kelly~\cite[\S3 and~\S4]{Kelly:Operads}):
$$F\circ G=\int^{m\in\Sigma}F(m)\otimes G^{\otimes m}=\coprod_{m\ge0} F(m)\otimes_{\Sigma_m}G^{\otimes m}.$$
\xexam

Recall that a category~$I$ is {\it sifted\/} if for all finite sets~$k$ the diagonal functor $I\to I^k$ is cofinal.
Filtered categories are sifted.
An example of a sifted category that is not filtered is given by the walking reflexive pair category,
consisting of two objects $0$~and~$1$ with two parallel arrows $f,g\colon0\to1$, and another arrow $h\colon1\to0$ such that $fh=gh=\id_1$.
Sifted colimits of this type are precisely reflexive coequalizers.
Any colimit can be expressed using reflexive coequalizers and coproducts, which explains why reflexive coequalizers appear constantly in constructions involving monoids and algebras over monoids.

\prop \label{prop--substitution.product}
The substitution product \refeq{substitution.product} is associative and unital.
Moreover, it is cocontinuous in the first variable and preserves sifted colimits in the second variable.
In particular, the substitution product is right closed, i.e., the functor~$-\circ G$ has a right adjoint for any~$G$.
\xprop

\pf
See \cite[Proposition~10.9 and Theorem~14.8]{GambinoJoyal:Operads}.
The bicategory of distributors used there is the opposite of the bicategory of finite sets, symmetric collections (with $\circ$ as the composition) and morphisms of collections.
\xpf

We emphasize that the substitution product does not preserve nonsifted colimits in the second variable,
for example, coproducts, because the functor $X\mapsto X^{\otimes k}$ in general does not preserve nonsifted colimits.
In particular, the substitution product is not left closed.
The substitution product is also not braided (in particular, not symmetric).
Note that the definition of the associator of~$\circ$ in the nonsymmetric case needs $\C$ to be \emph{symmetric} monoidal (see Muro~\cite[Remark~2.2]{Muro:Homotopy}).

\defi
The category $\ssOper := \ssOper_W(\C)$ of {\it $W$-colored (symmetric) operads in~$\C$\/} is the category of monoids in~$(\ssColl_W\C, \circ)$, i.e., $O \in \ssColl_W \C$ together with a unit map $1[1] \r O$ and a multiplication map $O \circ O \r O$ satisfying the associativity and unitality conditions.
For any set~$V$, the category of {\it $V$-colored (symmetric) modules\/} over a (symmetric) $W$-colored operad~$O$ is the category of left modules over~$O$ in~$\ssColl_{V,W}(\C)$.
It is denoted by $\Mod_O^V$.
Explicitly, its objects are given by $M \in \ssColl_{V,W}(\C)$ together with a map $O \circ M \r M$ subject to the standard associativity and unitality requirements.
For $V = \emptyset$ and $V=W$, we speak of \emph{$O$-algebras} and \emph{$O$-modules}, respectively, and denote them by $\Alg_O$ and $\Mod_O$.
Note that any $O$-algebra is naturally an $O$-module whose non-zero degrees are $\emptyset$.
\xdefi

The following result describes the categorical properties of colored modules over colored operads.

\theo \label{theo--O.Mod}
Suppose $(\C, \t)$ is a symmetric monoidal category that is enriched over a symmetric monoidal category $\V$.
Fix two sets $V$~and~$W$, and a $W$-colored (symmetric) operad $O$ in~$\C$.
\begin{enumerate}[(i)]
\item
\label{item--limits.O.Mod}
If $\C$ is complete, then so is $\Mod^V_O$ and the forgetful functor $U\colon\Mod_O^V \r \ssColl_{V,W}$ creates limits.

\item
\label{item--sifted.colimits.O.Mod}
If $\C$ admits sifted colimits (respectively, filtered colimits or reflexive coequalizers), which are preserved in each variable by the monoidal product in~$\C$,
then $\Mod^V_O$ admits sifted colimits, which are created by~$U$.

\item
\label{item--colimits.O.Mod}
If $\C$ admits reflexive coequalizers, which are preserved in each variable by the monoidal product in~$\C$,
then $\Mod^V_O$ is cocomplete.

\item
\label{item--locally.presentable.O.Mod}
If $\C$ is locally presentable and $\otimes$ preserves filtered colimits in each variable,
then $\Mod^V_O$ is locally presentable.

\item
\label{item--pullback.pushforward.O.Mod}
Suppose $f\colon O\to P$ is a morphism of $W$-colored (symmetric) operads in~$\C$.
If $\C$ admits reflexive coequalizers that are preserved in each variable by the monoidal product in~$\C$,
then the pullback functor $f^*\colon\Mod^V_P\to\Mod^V_O$ admits a left adjoint~$f_*$.
\end{enumerate}
\xtheo

\pf
Via \refpr{substitution.product}, these statements are reduced to similar statements about modules in (nonsymmetric, nonbraided) monoidal categories.
\refit{limits.O.Mod}, \refit{locally.presentable.O.Mod}, and \refit{pullback.pushforward.O.Mod} are then special cases of \cite[Theorem~3.4.1]{BarrWells:Toposes}, \cite[Theorem~5.5.9]{Borceux:2}, and \cite[Corollary~1]{Linton:Coequalizers}, respectively.

\refit{sifted.colimits.O.Mod}~\cite[Proposition~4.3.2]{Borceux:2} implies that $\Mod_O^V$ has sifted colimits, which are preserved by~$U$.
Reflection of sifted colimits by~$U$ is then implied by~\cite[Proposition~2.9.7]{Borceux:1} applied
to the opposite functor~$U^\op\colon(\Mod_O^V)^\op \to (\Mod_O^V)^\op$.
The cases of filtered colimits and reflexive coequalizers are treated identically.

\refit{colimits.O.Mod}~By \refit{sifted.colimits.O.Mod}, $\Mod_O^V$ admits reflexive coequalizers, which are created by~$U$.
Now apply~\cite[Corollary~2]{Linton:Coequalizers}, which in our case says that $\Mod_O^V$ has small colimits
if it has reflexive coequalizers and $\ssColl_{V,W}(\C)$ has small coproducts.
\xpf

\section{The enveloping operad}
\label{sect--enveloping.operad}

The enveloping operad (see, for example, \cite[Propositions~1.5]{BergerMoerdijk:Derived}, \cite[Proposition~5.4]{BergerMoerdijk:Axiomatic}) turns a module or algebra over an operad back into an operad.
This is used to relate properties of operadic algebras to those of operads, for example, pushouts (\refpr{filtration.Harper}) and transports along weak monoidal Quillen adjunctions (see \refth{transport}\refit{adjunction.Alg} and its proof).
We continue using the notation of \refsect{colored.collections}.

\defi
The category $\Pairs$ consists of pairs~$(O,A)$, where $O \in \ssOper_W$ is a (symmetric) $W$-colored operad in $\C$ and $A \in \ssColl_W$ is an $O$-module,
and a morphism of pairs~$(O,A)\to(P,B)$ is a morphism $f\colon O\to P$ of operads
together with a morphism $g\colon A\to f^*B$ of $O$-modules, where $f^*$ is the restriction functor from $P$- to $O$-modules.
\xdefi

\lemm
\label{lemm--enveloping.monoid}
There are adjunctions
$$\ssColl_W \mathrel{\mathop\rightleftarrows_{U}^{1[1] \x \id}} \Pairs \mathrel{\mathop\rightleftarrows_{\id \x U}^{\Env}} \ssOper_W.\eqlabel{adjunction.Pairs.Mon}$$
The functor $\id \x U$ sends an operad $O$ to $(O, U(O))$, where $U(O)$ is regarded as an $O$-module in the obvious way.
The functor $1[1] \x \id$ sends $X$ to $(1[1], X)$, where $1[1]$ is the initial operad.
The functor $U$ at the left sends $(O, M)$ to $U(M)$, i.e., it forgets the $O$-module structure on $M$.
The functor $\Env$ is called the \emph{enveloping operad}.
It satisfies $\Env(1[1], X) = \Free(X)$, where $\Free : \ssColl_W \rightleftarrows \ssOper_W : U$ is the free-forgetful adjunction.
\xlemm

\pf
The left adjunction holds since
$$\Pairs((1[1], X), (O, M)) = \ssColl_W(X, \eta^* M) = \ssColl_W(X, U(M)).$$
Here $\eta\colon 1[1] \r O$ is the unit of $O$, which is the unique morphism of operads $1[1] \r O$.
The right adjunction is a special case of \refth{O.Mod}\refit{pullback.pushforward.O.Mod} since $\Pairs$ are algebras over an operad similar to the operad of operads (\refsect{operad.operads}).
The last statement follows from the two adjunctions.
\xpf

\prop
\label{prop--enveloping.operad}
Fix a (symmetric) operad $O$ and consider the functor $\Env(O, -)\colon \Mod_O\r \ssOper_W$.
(We also apply this functor to $O$-algebras.)
\begin{enumerate}[(i)]
\item
\label{item--O.A.initial}
The enveloping operad of the initial $O$-algebra is given by $\Env(O, O \circ \emptyset) = O$.

\item
\label{item--O.A.connected.colimits}
The enveloping operad functor $\Env(O, -)$ preserves connected colimits of $O$-algebras, in particular, transfinite compositions.

\item
\label{item--O.A.pushouts}
Given a map $x\colon X \r X'$ in $\ssColl_W$, an $O$-module $A$, and a map $X \r U(A)$ in $\ssColl_W$, we form the pushout square in $\Mod_O$,
$$\xymatrix{
O \circ X \ar[d]^{O \circ x} \ar[r]^f & A \ar[d]^a \cr
O \circ X' \ar[r] & A'.\cr
}\eqlabel{O.X.A}$$
Then the following diagram is cocartesian in $\ssOper_W$,
where the top horizontal map is $$\Free(X) \buildrel\hbox{\ref{lemm--enveloping.monoid}}\over= \Env(1[1], X) \To7{\Env(\eta, f)} \Env(O, A):$$
$$\xymatrix{
\Free (X) \ar[r] \ar[d]^{\Free(x)}
&
\Free (U(A)) \ar[d]^{\Free(u)} \ar[r]^{\rm counit} &
\Env(O,A) \ar[d] \cr
\Free (X') \ar[r] &
\Free (U(A) \sqcup_X X') \ar[r] &
\Env(O,A').\cr
}
\eqlabel{free.operad}
$$
\item
\label{item--O.A.undercategory}
For any $A \in \Alg_O$, there is an equivalence of categories with the undercategory of $A$ in $\Alg_O$:
$$\Alg_{\Env(O,A)} = A \downarrow \Alg_O.$$
In particular, $\Env(O, A)_0 = A$.
\end{enumerate}
\xprop

\pf\leavevmode

\begin{enumerate}[(i)]
\item 
For any operad $T$, we have by adjunction
$$\ssOper_W(\Env(O, O\t \emptyset), T) = \{(f \in \ssOper(O, T), g\colon O \circ \emptyset \r f^* U(T) \in \Alg_O)\}.$$
As $O \circ \emptyset$ is initial in $\Alg_O$, $g$ is unique, so that this set of homorphisms is isomorphic to $\ssOper_W(O, T)$.
Hence, our claim.

\item 
For a connected index category $I$, $O$ is the colimit of the constant diagram $i \mapsto O$.
Therefore, $$(O, \colim A_i) = \colim (O, A_i).$$
Now apply the cocontinuity of the enveloping operad functor $\Pairs \r \ssOper$.

\item 
By \refle{enveloping.monoid}, the diagram \refeq{free.operad} is obtained by applying $\Env$ to the following diagram of pairs, which is easily seen to be cocartesian.
We conclude using that $\Env$ preserves all colimits, in particular pushouts.
$$\xymatrix{
(1[1], X) \ar[r]^{(1[1], f)} \ar[d]^{(1[1], x)} &
(1[1], U(A)) \ar[d] \ar[r]^{(\eta, \id)} &
(O, A) \ar[d] \\
(1[1], X') \ar[r] &
(1[1], U(A) \sqcup_X X')) \ar[r] &
(O, A' = A \sqcup_{O \circ X} O \circ X').
}
$$

\item 
Since the monoidal product in $\ssColl_W \C$ is right closed, an $\Env(O,A)$-module structure on some $X \in \ssColl_W$ is the same as a morphism of operads~$\Env(O,A)\to\End(X)$, where $\End(X) := \IHom(X, X) \in \ssOper_W$ is the endomorphism operad.
The adjunction \refeq{adjunction.Pairs.Mon} tells us that morphisms $\Env(O,A)\to\End(X)$ correspond to morphisms of pairs $(O,A)\to(\End(X),U (\End(X)))$.
This is the same as an $O$-module structure on $X$ and a map $A \r \End(X)$ of $O$-modules, where $\End(X)$ is regarded as an $O$-module via the chosen $O$-module structure on $X$.
Giving $A \r \End(X)$ is the same as $A = A \circ X \r X$.
The last equality uses that $A$ is an algebra, i.e., concentrated in degree $0$.

The second claim holds since $\Env(O, A)_0 = \Env(O,A) \circ \emptyset$ is the initial $\Env(O,A)$-module, which by the previous step is $A$.
\end{enumerate}
\xpf

\section{Admissibility of operads}
\label{sect--operad.admissible}

The following definition of admissibility of operads is standard (see, for example, \cite[\S2]{BergerMoerdijk:Derived}).

\defi \label{defi--operad.admissible}
A $W$-colored (symmetric) operad~$O$ in a symmetric monoidal model category~$\C$ is {\it admissible\/} if the product model structure on $\C^W$ transfers to $\Alg_O$ via the forgetful functor
$$\C^W \leftarrow \Alg_O : U,$$
i.e., if the classes $\we_{\Alg_O} = U^{-1}(\we_{\C^W})$ of weak equivalences and $\fib_{\Alg_O} = U^{-1}(\fib_{\C^W})$ of fibrations define a model category structure on $\Alg_O$.
Moreover, $O$ is \emph{strongly admissible} if it is admissible and if in addition $U$ preserves cofibrations with cofibrant source, i.e., for a cofibration $a\colon A \r A'$ of $O$-algebras, $U(a)$ is a cofibration and $U(A)$ is cofibrant in $\C^W$.
\xdefi

The admissibility of symmetric operads is a central problem in homotopical algebra.
It was addressed by Berger and Moerdijk \cite[Theorem~3.2]{BergerMoerdijk:Axiomatic} using the path object argument.
Their theorem requires the existence of a symmetric monoidal fibrant replacement functor and the monoidal unit to be cofibrant.
A well-known result due to Lewis \cite[Theorem~1.1]{Lewis:Is} precludes the existence of such data for a stable monoidal model category of spectra.
The conditions of their theorem were weakened by Kro \cite[Corollary~2.7]{Kro:Model}, whose version does not require the monoidal unit to be cofibrant.
Previously, Spitzweck had shown the existence of a semi-model structure for operads whose underlying symmetric sequence is projectively cofibrant
(which roughly means that $\Sigma_n$ acts freely on $O_n$) \cite[Theorem~4.7]{Spitzweck:Operads}.
This covers the Barratt--Eccles operad, for example, which satisfies $O_n=\E\Sigma_n$,
but excludes, say, the commutative operad $\Comm$, which is given by $\Comm_n=1$, the monoidal unit.
This is one of the most important examples of a symmetric operad, since its algebras are commutative monoid objects.
The admissibility of $\Comm$, i.e., the model structure on commutative monoid objects in~$\C$,
was established by Harper \cite[Proposition~4.20]{Harper:Symmetric}
and Lurie \cite[Proposition~4.5.4.6]{Lurie:HA} if $\C$ is freely powered.
Their proofs actually only use the weaker condition that the map $f^{\pp n}_{\Sigma_n}$ is an acyclic cofibration whenever $f$~is.
This property was later called the \emph{commutative monoid axiom} by White,
who also suggested a weakening similar to the one discussed in \refre{explain.assumptions} \cite[Theorem~3.2, Remark~3.3]{White:Model}.

The admissibility of arbitrary operads was also shown by Harper under the hypothesis that all objects in $\Sigma_n \C$ are \emph{projectively} cofibrant.
Again this is much stronger than being symmetric h-monoidal (see \csy{\refre{power}, \refchap{examples}}).
Subsequently to the present paper, White and Yau reproduced the admissibility of arbitrary operads under the condition that $X \t_{\Sigma_n} f^{\pp n}$ is an (acyclic) cofibration when $f$ is \cite[Theorem~6.1.1]{WhiteYau:Bousfield}.
This is a stronger assumption than symmetric h-monoidality, and is inapplicable to various flavors of spectra (e.g., symmetric, orthogonal, etc.)\ and other constructions used in stable homotopy theory, e.g., L-spaces.

For \emph{nonsymmetric operads}, the situation is quite a bit simpler, since no modding out by $\Sigma_n$ occurs in the definition of the circle product on nonsymmetric sequences.
Muro has shown the admissibility of all nonsymmetric operads under assumptions on $\C$ \cite[Theorem~1.2]{Muro:Homotopy}, \cite{Muro:Corrections},
which by \csy{\refle{i.monoidal.monoid.axiom}} are very closely related to the nonsymmetric part of \refth{O.Alg} (see \refre{explain.assumptions}).
Another type of admissibility result is due to Batanin and Berger who showed that so-called tame polynomial monads (or the colored operads associated to them) are admissible if $\C$ satisfies the monoid axiom \cite[Theorem 8.1]{BataninBerger:Homotopy}.

A technical key part in all proofs below is the analysis of pushouts of free $O$-algebra maps and free operad maps.
We will start with pushouts of operads and then deduce the pushouts of algebras from this.
The following description of pushouts of free (symmetric) operads is due to Spitzweck \cite[Proposition~3.5]{Spitzweck:Operads}
and, in the slightly different formulation given below, to Berger and Moerdijk \cite[Lemma~3.1]{BergerMoerdijk:Derived}, \cite[\S5.11]{BergerMoerdijk:Axiomatic}.

The description of such pushouts is based on the groupoid $\ssTree_W$ of \emph{$W$-colored (symmetric) marked trees}.
These are finite planar trees whose edges are labeled with colors $w \in W$.
The root vertex has a half-open (i.e., having only one boundary vertex) outgoing edge without called the \emph{root edge}.
It also has a (finite) number of vertices having half-open ingoing edges called the \emph{input edges}.
Any edge that is not a root edge nor an input edge is called an internal edge.
Their boundary consists of two vertices.
Moreover, a (finite) number of vertices of the tree is \emph{marked}, the others are not marked.
The markings is required to be such that every internal edge has at least one marked vertex at its boundary.
\emph{Automorphisms} of symmetric trees are isomorphisms of trees that do not respect the planar structure,
but do respect the markings, the colors of the edges, and send input edges to input edges.
Automorphisms of nonsymmetric trees are only identity morphisms.
For a vertex $r$ in a tree, the \emph{valency} $\val(r) \in \ssSeq_W$ is given by $(s, w)$, where the multisource $s\colon I \r W$ is given by the set~$I$ of the incoming edges of~$r$,
ordered according to the planar structure (which is only needed to make this notion unambiguous) and their corresponding colors, and target~$w$ given by the color of the outgoing edge.
In the same vein, the valency $\val(T)$ of the tree is given by the colors of the input edges and root edge.
The subgroupoid of trees with $k$ marked vertices and valency $(s, w) \in \ssSeq_W$ is denoted by $\ssTree_{s, w}^{(k)}$.

Using the notation of \refpr{filtration.operads}, the intuitive meaning of these notions is that a tree~$T$ with valency~$(s, w)$ stands for an operation in~$O'$ with inputs given by the multisource~$s$ and target~$w$.
Such operations are nested applications of the more elementary operations given by vertices.
If $T$ contains no marked vertices, i.e., $k = 0$, then $T$ is just a corolla consisting of a root edge and finitely many input edges, corresponding to the operations that are present in $O$.
More generally, for $k \ge 0$, $k$ operations coming from $\Free(X)$ have been identified by their image in $\Free(X')$.

\prop (Spitzweck, Berger--Moerdijk)
\label{prop--filtration.operads}
Let $\C$ be a symmetric monoidal model category.
For any map $x\colon X \r X'$ in $\ssColl_W$ and any pushout diagram in $\ssOper_W$,
$$\xymatrix{
\Free(X) \ar[d]_{\Free(x)} \ar[r] & O \ar[d]^o\cr
\Free(X') \ar[r] & O',\cr
}\eqlabel{pushout.operads}$$
the map $U(o)_{s,w} \in \Ax_s \C$ is the transfinite composition of maps $O^{(k)}_{s, w} \r O^{(k+1)}_{s, w}$, for $k \ge 0$, which arise as the following pushouts in $\Ax_s \C$:
$$\xymatrix{
\coprod_{T} \Ax_s \cdot_{\Aut T} x^*(T) \ar[r] \ar[d]_{\coprod_T \Ax_s \cdot_{\Aut T} \epsilon(T)} & O^{(k)}_{s,w} \ar[d]\cr
\coprod_{T} \Ax_s \cdot_{\Aut T} x(T) \ar[r] & O^{(k+1)}_{s,w}.\cr
}\eqlabel{filtration.operad}$$
The coproducts run over all isomorphism classes of $\ssTree^{(k)}_{s,w}$ as defined above.
For such a tree $T$, the map $\epsilon(T)\colon x^*(T) \r x(T)$ is inductively defined as
$$\epsilon(T) := \epsilon(r(T)) \pp \underbrace{\bigpp_i \epsilon(T_i)^{\pp t_i}}_{=:\epsilon'(T)},$$
where $\epsilon(r(T)) \in \Ax_{\val(r(T))} \C$ is defined as
$$\epsilon(r(T)):=\cases{x_{\val(r(T))},&if $r(T)$ is marked;\cr
(\eta_O)_{\val(r(T))},&if $r(T)$ is not marked,\cr}\eqlabel{epsilon}$$
where
$\eta_O\colon 1[1] \r U(O)$ is the unit map of $O$ and $\val(r(T))$ is the valency of the root $r(T)$ of $T$.
Isomorphic subtrees (with markings, colors, and input edges induced from $T$) of the root are grouped together and denoted by $T_i$, $1 \le i \le k$.
The number of subtrees isomorphic to $T_i$ is denoted by~$t_i$, so that $\sum_{i=1}^k t_i$ equals the cardinality of the multisource of $r(T)$.
The group
$$\Aut (T) = \prod_{i=1}^k \Aut(T_i)^{t_i} \rtimes \prod_{i=1}^k \Sigma_{t_i}$$
acts on $\epsilon(r(T))$ via the quotient $\prod \Sigma_{t_i}$ and in the natural way on $\epsilon'(T) \in (\prod \Aut(T_i)^{t_i}) \C$.
\xprop

\pf
This is exactly the statement of Berger and Moerdijk cited above, if we replace $\epsilon(T)$ by $\epsilon_u(T)$, which is defined as above, except that $\epsilon(r(T)) := u_{\val(r(T))}$ if the vertex $r(T)$ is marked, where $u\colon U(O) \r U(O) \sqcup_X X'$ is the pushout of $x$.
We conclude using the pushout square $\Ax_s \cdot_{\Aut T} \epsilon(T) \r \Ax_s \cdot_{\Aut T} \epsilon_u(T)$ and \csy{\refpr{pushout.product.pushout.morphisms}}.
\xpf

\refpr{filtration.operads} has the following model-categorical consequence due to Spitzweck \cite[Lemma~3.6]{Spitzweck:Operads} and, in the form below, to Berger--Moerdijk \cite[Proposition~5.1]{BergerMoerdijk:Axiomatic}.
We will show in \refle{cofibrancy.preserved}\refit{forget.Oper} that $U(\eta_O)$ is a cofibration for any cofibrant operad $O$, so the corollary is applicable to such pushouts.
This will be important in the study of strong admissibility.
Recall that $\ssColl_W(\C)$ is equipped with the projective model structure.
Unless the contrary is explicitly stated, all cofibrations in categories of the form $G\C$, for a finite group $G$, are understood as \emph{projective} cofibrations.
(The distinction between injective and projective model structures only matters in the symmetric case, for the category of nonsymmetric collections $\nsColl_W(\C)$ is just a product of copies of $\C$.)

\coro
\label{coro--filtration.operads.model}
In the situation of \refpr{filtration.operads}, suppose that $U(\eta_O)$ is a cofibration in $\ssColl_W$.
Also suppose that $x$ is a cofibration in $\ssColl_W$.
Then the vertical maps in \refeq{filtration.operad} are cofibrations in $\Ax_s \C$.
Therefore, $U(o)$ is also a cofibration in $\ssColl_W$.
\xcoro

The following description of pushouts of free $O$-algebras is due to
Fresse \cite[Proposition~18.2.11]{Fresse:Modules}, Elmendorf and Mandell~\cite[\S12]{ElmendorfMandell:Rings},
and Harper \cite[Proposition~7.12]{Harper:Monoidal}.

\prop
\label{prop--filtration.Harper}
Let $\C$ be a symmetric monoidal model category and $O$ a (symmetric) operad.
Let $$\xymatrix{O \circ X \ar[d]_{O \circ x} \ar[r] & A \ar[d]^a\cr
O \circ X'\ar[r] & A'\cr}
\eqlabel{pushout.O.Alg}$$
be a pushout diagram of $O$-algebras, where $x\colon X \r X'$ is a map in $\C^W$.
For any color $w \in W$, the map $U(a)_w \in \C$ lies in the weak saturation of morphisms of the form
$$\Env(O, A)_{s, w} \t_{\Ax_s} \bigpp_{r \in W} x_r^{\pp s^{-1}(r)},\quad s\colon I \r W \in \ssSeq_W^\x,\quad I \ne \emptyset.\eqlabel{U.a}$$
(The pushout product is finite, since $I$ is a finite set.)
For example, if $W$ consists of a single color and we consider symmetric operads, $U(a)$ lies in
$$\ws{\{\Env(O, A)_n \t_{\Sigma_n} x^{\pp n}, \ n \ge 1 \}}.$$
\xprop

\pf
By \refpr{enveloping.operad}\refit{O.A.undercategory}, the map $U(a)_w$ is the level $(\emptyset, w)$ of $\Env(O, A) \r \Env(O, A')$,
which by the pushout diagram \refeq{free.operad} and description of pushouts in \refpr{filtration.operads} is a transfinite composition of pushouts of the maps \refeq{filtration.operad} (where the $O$ there is now $\Env(O,A)$).
The map $x$ is concentrated in degree 0, so the only trees $T$ such that the map $\epsilon(T)$ defined in \refeq{epsilon} is not an isomorphism are the trees (with valence $(\emptyset, w)$) whose marked vertices have valency 0, i.e., are stumps.
Since any internal edge has at least one marked vertex, the only such trees $T$ are corollas whose root is not marked and has valence $(t\colon I \r W, w)$ and whose leaves are marked.
We get $\epsilon(T) = \Env(O,A)_{t,w} \t \bigpp_{i \in I} x_{t(i)}$ and $\Aut (T)=\Ax_t$.
Hence, the left-hand vertical map in \refeq{filtration.operad} agrees with~\refeq{U.a}.
\xpf

The next result identifies (symmetric) h-monoidality as the key condition for admissibility of all (symmetric) operads.
We emphasize that symmetric h-monoidality requirement is stable under weak saturation, transfer of model structures, and monoidal left Bousfield localization (see \csy{Theorems \ref{theo--symmetric.weakly.saturated}, \ref{theo--transfer.symmetric.monoidal}, and \ref{theo--symmetric.monoidal.localization}} for the precise statements).
Basic examples of symmetric h-monoidal model categories include simplicial sets, simplicial presheaves, topological spaces,
chain complexes of rational vector spaces, and symmetric spectra (see \csy{\refchap{examples}}).
Chain complexes of abelian groups are not symmetric h-monoidal and, in fact, the commutative operad is provably not admissible in chain complexes of abelian groups.
Recall the definitions of the terms below from \refde{summary.symmetricity}.

\lemm
\label{lemm--technical}
Suppose $\C$ is a model category.
\begin{enumerate}[(i)]
\item
If the class of weak equivalences in~$\C$ is stable under colimits of chains, then the same is true for the class of h-cofibrations,
and, if $\C$ is left proper, for acyclic h-cofibrations.
\item
If a symmetric monoidal category~$\C$ is left proper, strongly admissibly generated, and satisfies the acyclic part of symmetric h-monoidality,
then any transfinite composition of pushouts of maps of the form $Y \t_{\Ax_s} x^{\pp s}$ is an acyclic h-cofibration,
where $Y \in \Ax_s \C$ is arbitrary and $x$ is a finite family of acyclic cofibrations.
\end{enumerate}
\xlemm

\pf
The first part is \csy{\refle{i.cofibrations}\refit{i.cofibrations.weakly.sat}}
and the second part is \csy{\refpr{strongly.admissibly.generated}} (first paragraph in the proof).
For the convenience of the reader, we sketch the argument in the latter case.
The maps in question lie in the class \refeq{cellular.algebraic} and are acyclic h-cofibrations by acyclic symmetric h-monoidality.
Hence, by assumption, the (co)domains of generating cofibrations of~$\C$ are compact with respect to the cellular closure of these maps.
Analogously to \csy{\refle{sequential}\refit{we.chain.colimits}}, a transfinite composition~$f_\infty$ of weak equivalences $f_i$ (in our case, cobase changes of the above maps)
is a weak equivalence,
provided that (co)domains of the generating cofibrations of~$\C$ are compact relative to the class spanned by the acyclic cofibrations and the maps in the transfinite chain.
Similarly, $f_\infty$ is an h-cofibration provided that the $f_i$ and the maps in the chain are h-cofibrations.
\xpf

\theo \label{theo--O.Alg}
Suppose $\C$ is a symmetric monoidal model category and $W$ is a set.
Furthermore, suppose that one of the following is satisfied:
\begin{enumerate}[(a)]
\item
\label{item--combinatorial.case}
$\C$ is combinatorial and weak equivalences are closed under transfinite compositions, or
\item
\label{item--admissibly.generated.case}
$\C$~is admissibly generated and tractable.
\end{enumerate}
If $\C$ is (symmetric) h-monoidal (the acyclic part is sufficient), then any $W$-colored (symmetric) operad~$O$ in~$\C$ is admissible.
\xtheo

\pf
We apply \cite[Theorem~11.3.2]{Hirschhorn:Model} to the adjunction $O \circ - : \C^W \rightleftarrows \Alg_O : U$.
By \refth{O.Mod}, $U$ preserves sifted colimits and $\Alg_O$ is complete and cocomplete.

We now show that transfinite compositions of the images under~$U$ of cobase changes of elements in~$F(J)$ are weak equivalences in~$\C^W$.
Consider a cocartesian diagram of $O$-algebras as in \refeq{pushout.O.Alg},
where $x\colon X\to X'$ is generating acyclic cofibration in~$\C^W$, hence also an acyclic (symmetric) h-cofibration.
By \refpr{filtration.Harper}, the morphism~$U(a)$ is the (countable) transfinite composition of cobase changes of morphisms
$$\Env(O, A)_{s, w} \t_{\Ax_s} \bigpp_{r \in W} x_r^{\pp s^{-1}(r)}, \qquad s\colon I \r W \in \ssSeq_W^\x.\eqlabel{proof.Harper}$$
Here, $\Env$ is the enveloping operad (\refle{enveloping.monoid}) and $\Ax_s$ is the group of automorphisms of the multisource~$s$,
which is trivial for nonsymmetric operads, and as in \refeq{Ax.s} for symmetric operads.
Each of the above morphisms is a couniversal weak equivalence or,
equivalently \cite[Lemmas 1.6 and~1.8]{BataninBerger:Homotopy}, an acyclic h-cofibration
since $x$ is an acyclic (symmetric) h-cofibration, i.e., each $x_r$ is one.
Their transfinite composition is again a couniversal weak equivalence by \refle{technical}.

We finally show that $F(I)$ and $F(J)$ permit the small object argument \cite[Definition~10.5.15]{Hirschhorn:Model}.
If $\C$ is combinatorial, this is tautological since all objects are small.
Suppose now that $\C$ is admissibly generated and tractable.
By \refde{summary.symmetricity}, all cofibrant objects, in particular the (co)domains of $I$ are small relative to $\cell(-)$ applied to the maps in \refeq{proof.Harper},
where $x$ is a cofibration.
Therefore, they are small relative to $U(\cell (O \circ I))$.
By adjunction, the (co)domains of $O \circ I$ are therefore small relative to $\cell(O \circ I)$.
Again using the tractability, the same argument shows that $O \circ J$ is small relative to $\cell(O \circ I)$, a fortiori relative to $\cell(O \circ J)$.
\xpf

\rema
\label{rema--explain.assumptions}
\label{rema--module.context}
The proof also shows the following statement: suppose $\C$ is a symmetric monoidal category, $\C'$ is a combinatorial (more generally,
admissibly generated) and such that $\C'$ is a commutative $\C$-algebra.
Finally, suppose that for a finite family of generating cofibrations $x_{r_1}, \ldots, x_{r_k}$ in $\C'$, and $n_1, \ldots, n_k \ge 1$, any object $E \in (\prod_{j=1}^k \Ax_{n_j}) \C$, the map
$$E \t_{\prod_j \Ax_{n_j}} \bigpp_j x_{r_j}^{\pp n_j}\eqlabel{explain.assumptions}$$
lies in a class whose saturation under transfinite composition and pushouts consists of weak equivalences (in~$\C'$).
Then any $W$-colored symmetric operad $O$ in $\C$ is admissible, i.e., the $O$-algebras in $\C'$ carry a transferred model structure.
Since the differences are purely grammatical, we omit the proof of this assertion.

The same statement holds for nonsymmetric operads after dropping $\prod \Ax_{n_j}$ in \refeq{explain.assumptions}.
If, in addition, the monoidal product of $\C'$ turns $\C'$ into a monoidal model category it can be further simplified to requiring the above condition only for the maps $E \t x$, where $E \in \C$ and $x \in \C'$ is a generating acyclic cofibration.
This is exactly the monoid axiom \cite[Definition~3.3]{SchwedeShipley:Algebras}, so the above proof reproduces one of Muro's aforementioned admissibility result of nonsymmetric operads
\cite[Theorem~1.2]{Muro:Homotopy} and~\cite{Muro:Corrections}.

In particular, the nonacyclic part of (symmetric) h-monoidality is not necessary for the admissibility statement.
We mention the nonacyclic part in the definition of (symmetric) h-monoidality,
since the combination of the acyclic and the nonacyclic part of (symmetric) h-monoidality is easier to localize.
Also, for concrete model categories, it is usually easier to establish both properties simultaneously.
For the same reason, we have separated the saturation with respect to transfinite compositions
and the one with respect to pushouts (governed by (symmetric) h-monoidality).
See \csy{\refth{symmetric.monoidal.localization}\refit{localization.symmetric.i.monoidal}, \refchap{examples}} and the remarks at the end of \refsect{symmetricity}.
\xrema

\section{Strong admissibility of operads}
\label{sect--strong.admissibility}

In addition to the admissibility of operads, it is in practice desirable to know when the forgetful functor
$$\C^W \gets \Alg_O: U$$ preserves cofibrant objects or even cofibrations with cofibrant source, i.e., when $O$ is strongly admissible,
as defined in \refde{operad.admissible}.
We present two results in this direction: \refpr{cofibrant.strongly.admissible} is a result for levelwise projectively cofibrant operads.
It works in any symmetric monoidal model category.
\refth{O.strongly.admissible} is a much more flexible criterion for levelwise \emph{injectively} cofibrant operads.
Here, the additional key condition is the symmetroidality of $\C$.
More precisely, we use the following conditions on $O$ (well-pointedness, also known as $\Sigma$-cofibrancy was considered in \cite{BergerMoerdijk:Axiomatic, BergerMoerdijk:Derived}):

\defi
\label{defi--well.pointed}
A (symmetric) colored operad $O$ is \emph{well-pointed} (respectively, \emph{weakly well-pointed})
if the unit map $U(\eta_O)$ is a projective (respectively, injective) cofibration in $\ssColl_W(\C)$.
\xdefi

Strong admissibility does not seem to have been studied before as an independent notion.
See Mandell \cite[Lemma~13.6]{Mandell:Einfty}, Shipley \cite[Proposition~4.1]{Shipley:Convenient}, and Harper and Hess \cite[Theorem~5.18]{HarperHess:Homotopy}
for strong admissibility statements for operads in chain complexes, simplicial symmetric spectra, and arbitrary monoidal model categories, though.

The following preparatory lemma captures the preservation of cofibrant objects under various forgetful functors.
We do not claim originality for this lemma, for example, Part~\refit{Env.preserves.cofibrations} is similar to \cite[Proposition~2.3]{BergerMoerdijk:Derived}.

\lemm \label{lemm--cofibrancy.preserved}
With $\C$~and~$W$ as before, the following claims hold:
\begin{enumerate}[(i)]
\item \label{item--forget.Oper}
Let $f\colon O \r O'$ be a cofibration in $\ssOper_W$ such that $O$ is well-pointed, i.e., $U(\eta_O)$ is a cofibration in $\ssColl_W$.
Then $U(f)$ is a cofibration in $\ssColl_W$.
In particular,
\begin{enumerate}[(1)]
\item \label{item--special.case.1}
any cofibrant operad $O$ is well-pointed,
in other words, the levels $O_{s,w}$ are cofibrant in $\Ax_s^\proj \C$ for all $s\colon I \r W$ if $\sharp I \ne 1$ or if $\sharp I = 1$ and $s(*) \ne w$ and the unit map $1 \r O_{w,w}$ is a cofibration in $\C$ for all $w \in W$;
\item \label{item--special.case.2}
the forgetful functor $U$ sends cofibrations with cofibrant source to cofibrations;
\item \label{item--special.case.3}
if the unit $1 \in \C$ is cofibrant, $U$ also preserves cofibrant objects, i.e., the underlying (symmetric) sequence $U(O) \in \ssColl_W$ of any cofibrant operad $O$ is cofibrant.
\end{enumerate}
\item \label{item--Env.preserves.cofibrations}
For any (symmetric) operad $O$, the functor $\Alg_O \r \ssOper_W$, $A \mapsto \Env(O, A)$ preserves cofibrations.
For example, $O \r \Env(O, A)$ is a cofibration for any cofibrant $O$-algebra $A$.
\end{enumerate}
\xlemm

\pf
\refit{forget.Oper}~The map $f$ is a retract of a transfinite composition of pushouts of maps $\Free(x)$ as in \refeq{pushout.operads},
where $x$ is a cofibration in $\ssColl_W$ and, by assumption and cellular induction, $O$ is well-pointed. 
The functor $U$ commutes with retracts and transfinite compositions.
Cofibrations (in $\ssColl_W$) are stable under these two types of saturation.
Therefore the statement follows from \refco{filtration.operads.model}, using that $O$ is well-pointed. 

The remaining statements are special cases: \refit{special.case.1} follows by applying the general statement to $\eta_O\colon 1[1] \r\nobreak O$.
\refit{special.case.2} follows by combining the general statement and \refit{special.case.1}.
Finally, \refit{special.case.3} holds since $1[1]$ is the initial operad, whose underlying symmetric sequence is cofibrant in $\ssColl_W$ if and only if $1$ is cofibrant in $\C$.

\refit{Env.preserves.cofibrations}~The claim about $\Env(O, -)$ follows from \refpr{enveloping.operad}:
if $a$ is a pushout of a free $O$-algebra map $O \circ x$ on a cofibration $x \in \ssColl_W$ as in \refeq{pushout.O.Alg}, the map $\Env(O, a)\colon \Env(O, A) \r \Env(O, A')$ is the pushout of $\Free(x)$, which is a cofibration (in $\ssOper_W$).
For a transfinite composition of cofibrations of $O$-algebras, we use that both $U$ and $\Env(O, -)$ preserve filtered colimits.
By \refpr{enveloping.operad}\refit{O.A.initial}, the last statement is the special case $a\colon O_0 = O \circ \emptyset \r A$.
\xpf

The following result guarantees strong admissibility for those operads whose levels are projectively cofibrant (except for unit degrees, in which case the map from the monoidal unit to the level is required to be a cofibration).
By \cite[Theorem~4]{Spitzweck:Operads}, any cofibrant operad $O$ is admissible if $\C$ satisfies the monoid axiom, so it is strongly admissible in this case by the result below.

\prop \label{prop--cofibrant.strongly.admissible}
Suppose $\C$ is a symmetric monoidal model category.
Any admissible well-pointed (symmetric) operad $O \in \ssOper_W(\C)$ is strongly admissible.
For example, any admissible cofibrant operad is strongly admissible.
\xprop
\pf
Let $E := O \circ \emptyset = O_0$ be the initial $O$-algebra.
For any two cofibrations $E \r A \r B$ of $O$-algebras,
both maps $O \stackrel{\text{\ref{prop--enveloping.operad}\refit{O.A.initial}}}= \Env(O, E) \stackrel \epsilon \r \Env(O, A) \stackrel \alpha \r \Env(O, B)$ are cofibrations of operads by \refle{cofibrancy.preserved}\refit{Env.preserves.cofibrations}.
By assumption and \refle{cofibrancy.preserved}\refit{forget.Oper}, $U(\epsilon)$ is a cofibration in $\ssColl_W(\C)$.
Therefore, $\Env(O, A)$ is again well-pointed, so that $U(\alpha)$ is a cofibration for the same reason.
In particular, the zeroth level of $U(\alpha)$, which by \refpr{enveloping.operad}\refit{O.A.undercategory} is $U(A) \r U(B)$, is a cofibration in $\C^W$.
\xpf

The next theorem is a supplementary condition for strong admissibility of arbitrary symmetric operads.
Recall from \csy{\refchap{examples}} that rational chain complexes and symmetric spectra (with an appropriate stable positive model structure) are symmetroidal.
The latter statement also shows that under very mild conditions, any monoidal model category is Quillen equivalent to a symmetroidal model category.
Moreover, symmetroidality is stable under monoidal Bousfield localization and transfer (see \csy{\refth{transfer.symmetric.monoidal} and \refth{symmetric.monoidal.localization}} for the precise statements).
These results turn \refth{O.strongly.admissible} into a powerful tool ensuring strong admissibility of operads.

The following lemma is the key stepstone for strong admissibility.
In order to keep the exposition brief, we will again speak of ``(symmetric) operads'' in a symmetric monoidal category to simultaneously cover the case of symmetric and of nonsymmetric operads.
In the latter case, all the groups $\Ax_s$ and $\Aut T$ appearing below are trivial by definition.

\lemm
\label{lemm--preparation.strongly.admissible}
Let $\C$ be a symmetric monoidal model category.
Let $O$ be a (symmetric) $W$-colored operad and $A$ any cofibrant $O$-algebra.
For any $(s\colon I \r W,w) \in \ssSeq_W$, the levels of the unit map
$$(\eta_{\Env(O,A)})_{s,w}\colon 1[1]_{s,w} \r \Env(O,A)_{s,w}$$
in $\Ax_s \C$ are contained in $\ws{(\Y_O)_s}$, where $(\Y_O)_s$ is the smallest class of morphisms in $\Ax_s \C$ that contains all isomorphisms, the generating cofibrations of $\C$ (for $\sharp I = 0$ only), and finally contains
$$(\eta_O)_{s \sqcup t, w} \pp_{\Ax_t} x^{\pp t} := (\eta_O)_{s \sqcup t, w} \pp_{\Ax_t} \bigpp_r x_r^{\pp n_r}.\eqlabel{eta.O.x}$$
Here, $t\colon J \r W$ is any multisource and the multi-index $n$ is given by $n_r = \sharp t^{-1}(r)$ for $r \in W$, and $x=(x_r)$ is a finite family of generating cofibrations in $\C$.
(We use the convention that only the finitely many terms with $n_r \ne 0$ appear, unless $J = \emptyset$, in which case we interpret the above expression as $(\eta_O)_{s,w}$.)

In particular, for any cofibrant $O$-algebra $A$, the map $\emptyset \r U(A) \in \C^W$ is contained in
$$\ws{\cof_\C \cup \{ (\eta_O)_{t, w} \pp_{\Ax_t} x^{\pp t}, (t, w) \in \ssSeq_W \}}.$$
\xlemm

\pf
We prove this by cellular induction on $A$, using the properties of the enveloping operad established in \refpr{enveloping.operad}.
We will write $\varphi \colon G\C \r \C$ for any functor that forgets the action of some finite group~$G$, for example, $G=\Ax_s$.
For $A = O \circ \emptyset = O_0$, $O = \Env(O, O_0)$ is an isomorphism, so the claim is clear by assumption.
For a pushout of $O$-algebras as in \refeq{pushout.O.Alg}, where $A$ is cofibrant and $x$ is a cofibration, there is a pushout of operads
$$\xymatrix{
\Free(X) \ar[d]_{\Free(x)} \ar[r] & \Env(O,A) \ar[d]^o \\
\Free(X') \ar[r] & \Env(O,A').
}
\eqlabel{pushout.enveloping.operads}$$
We now use \refpr{filtration.operads}, including the notation.
We need to show
$$\Ax_s \cdot_{\Aut T} \epsilon(T) \in (\Y_O)_s.$$
By induction on the tree $T$, one sees that
$$\varphi (\epsilon(T)) = \bigpp_{r \in T} \varphi(\epsilon(r)),$$
where the pushout product runs over all vertices $r$ of $T$.
Recall that $f \pp g$ is an isomorphism for all maps $g$ whenever $f$ is an isomorphism.
Hence, it is enough to prove our claim for those trees $T$ such that none of the morphisms~$\epsilon(r)$ is an isomorphism.

If a vertex $r \in T$ is marked, then $\epsilon(r) = u_{\val(r)}$, where
$u\colon U(\Env(O,A)) \r U(\Env(O,A)) \sqcup_X X'$ is the pushout of $x$ along the map $X \r U(\Env(O,A))$ adjoint to the top horizontal map in \refeq{pushout.enveloping.operads}.
If $r$ is marked and has positive valency, i.e., $(s,w) := \val(r)$ with a multisource $s\colon I \r W$ of arity $\sharp I > 0$, then $u_{s,w}$, which is a pushout of $x_{s,w} = \id_{\emptyset}$, is an isomorphism.
Thus, we may assume that the marked vertices have valency~0, i.e., no incoming edges.
On the other hand, by definition of marked trees, any edge contains at most one nonmarked vertex.
Therefore, the only trees we need to consider are as follows.
\begin{enumerate}
\item The tree denoted by ${}^w+$ consisting of a single marked vertex with no incoming edge and the outgoing root edge colored by $w$.
\item The trees denoted by ${}^w-^{t+}_s$ consisting of a single nonmarked vertex that has a root edge of color $w$, some noninput edges whose other end is marked, and some input edges.
The valency of the input and noninput edges is denoted by $s$ and~$t$ respectively.
\end{enumerate}
Here is a picture of ${}^w+$ and of ${}^w-^{t+}_s$.
The different dashing styles indicate different colors, the two rightmost lower arrows are input edges, the top arrows are the root edges, $\bullet^+$ is a marked vertex, $\bullet^-$ is not marked.
$$
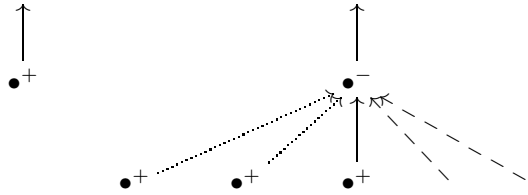
\xymatrix{
& & & & & \\
\bullet^+ \ar[u] & & & \bullet^- \ar[u] \\
& \bullet^+ \ar@{.>}[urr] & \bullet^+ \ar@{.>}[ur] & \bullet^+ \ar@{->}[u] & \ar@{-->}[ul] & \ar@{-->}[ull]
}$$
For $T = {}^w+$, we have $\Sigma_s = \Aut T = 1$ and $\epsilon(T) = x_w$, which is in $\Y_O$ being a cofibration.
For $T = {}^w-^{t+}_s$, we have $\Aut(T) = \Ax_{s} \x \Ax_{t}$, where $\Ax_s$ and $\Ax_t$ are defined in \refeq{Ax.s}.
In the example above, $\Sigma_{t} = \Sigma_2 \x \Sigma_1$ and $\Sigma_{s} = \Sigma_2$.
We group the noninput edges of $\bullet^-$ according to their color, say $n_i$ noninput edges of color~$t_i$.
Then
$$\Sigma_s \cdot_{\Aut T} \epsilon(T) = (\eta_{\Env(O,A)})_{s \sqcup t, w} \pp_{\prod \Sigma_{n_i}} \bigpp_{i} x_{t_i}^{\pp n_i},$$
which is in $\Y_O$ by the inductive hypothesis.
This finishes the pushout step.

The handling of retracts and transfinite compositions of cofibrant $O$-algebras is clear, noting that the functor $\Alg_O \r \ssColl_W$,
$A \mapsto U(\Env(O,A))$ preserves filtered colimits and retracts.

The claim concerning $U(A)$ is the restriction of the statement about the levels of $\Env(O,A)$ to degree~0.
\xpf

\theo \label{theo--O.strongly.admissible}
Suppose that $\C$ is a symmetric monoidal model category and $O$ is an admissible (symmetric) $W$-colored operad in $\C$.

In the nonsymmetric case, suppose that $(\eta_O)_{s,w} \pp - \colon \Ar(\C) \r \Ar(\C)$ preserves (acyclic) cofibrations.

In the symmetric case, suppose that $\C$ is symmetroidal (\refde{summary.symmetricity}) with respect to the class $\Y_O = ((\Y_O)_n)$ consisting of
$$(\Y_O)_n := \bigcup_{(s, w)} (\Y_O)_s,$$
where, as above, $s$ is such that $n_r = \sharp s^{-1}(r)$ (for $r \in W$), $w \in W$ is arbitrary, and $(\Y_O)_s$ is the class of morphisms in $\Ax_s \C$ defined in \refle{preparation.strongly.admissible}.

Then $O$ is strongly admissible.

For example, if $\C$ is symmetroidal (i.e., symmetroidal with respect to the injective cofibrations in $\Sigma_n \C$) and $O$ is \emph{weakly} well-pointed (i.e., $O_{s,w}$ is cofibrant in $\C$ for all nonunit degrees $(s,w)$ and $1 \r O_{w,w}$ is a cofibration in $\C$), then $O$ is strongly admissible.
\xtheo

\pf
It is enough to show that the maps in \refeq{U.a} are cofibrations in $\C^W$ for any cofibrant $O$-algebra $A$ and any cofibration $x$ in $\C^W$.

To show this in the symmetric case, by the symmetroidality condition on $\C$ and \csy{\refle{symmetroidal.Y.weakly.sat}}, which allows to weakly saturate the symmetroidality class, we have to show that the map
$$(\eta_{\Env(O,A})_{s,w} \colon 1[1]_{s,w} =\cases{1,&unit degrees;\cr\emptyset,&nonunit degrees.\cr}\To2{}\Env(O,A)_{s,w}$$
lies (levelwise) in $(\Y_O)_s$.
For unit degrees $(s,w) = (w,w)$, this guarantees that $\Env(O,A)_{w,w} \t x$ is a cofibration by \refle{monoidally.cofibrant.pushout.product}\refit{atB}.
This is exactly the content of \refle{preparation.strongly.admissible}.

In the nonsymmetric case, the argument is similar, but considerably easier since $\Ax_s$ is trivial: if the pushout product with $(\eta_O)_{s,w}$ preserves (acyclic) cofibrations, then so does the pushout product with the maps in \refeq{eta.O.x} and therefore also the pushout product with $(\eta_{\Env(O,A)})_{s,w}$.
Again, this implies that the maps in \refeq{U.a} are cofibrations in $\C^W$.

The last statement is a special case: let $\C$ be symmetroidal, i.e., symmetroidal with respect to $\Y_n := \cofib_{\Sigma_n^\inje \C}$.
Then $\Y_n \supseteq (\Y_O)_n$: indeed, the maps in \refeq{eta.O.x} are injective cofibrations by the symmetroidality of $\C$.
\xpf

The following corollary illustrates how to transfer the strong admissibility of operads.
Note that the symmetroidality of $\C$ does \emph{not} imply the symmetroidality of $\D$, i.e., the symmetroidality with respect to $\cofib_{\Sigma_n^\inje \D}$,
but only the symmetroidality with respect to $F(\cofib_{\Sigma_n^\inje \C})$ (see \csy{\refth{transfer.symmetric.monoidal}\refit{transfer.symmetroidal} and \refre{issue.symmetroidal}}).

\coro \label{coro--strongly.admissible.transfer}
Let $F: \C \leftrightarrows \D : G$ be a Quillen adjunction of symmetric monoidal model categories such that the model structure on $\D$ is transferred from $\C$ and such that $F$ is strong symmetric monoidal.
Suppose $\C$ is symmetroidal (only required in the symmetric case) and let $O$ be a weakly well-pointed (symmetric) operad in~$\C$.
We assume that the operad $P$ in $\D$ given by $P_{s,w} = F(O_{s,w})$ is admissible.
Then $P$ is strongly admissible.
\xcoro

\pf
The strong monoidality of $F$ gives the strong monoidality of the left adjoint in the adjunction
$$F: (\ssColl_W(\C), \circ) \rightleftarrows (\ssColl_W(\D), \circ): G.$$
The resulting adjunction of monoids, i.e., $W$-colored operads (see also \refeq{adjunction.Oper})
$$F^\ssOper: \ssOper_W(\C) \rightleftarrows \ssOper_W(\D): G$$
is therefore such that $U_\D F^\ssOper = F U_\C$, where $U_?: \sColl_W(?) \r ?$ are the forgetful functors.
Therefore, $P$ as defined above, is indeed an operad.

As in the proof of \refth{O.strongly.admissible}, we have to show that $\D$ is $\Y_P$-symmetroidal.
The generating cofibrations $y$ of $\D$ are of the form $y=F(x)$, $x \in \cof_\C$.
The (levels of) $U(\eta_P)$ are of the form $FU(\eta_O)$.
Finally, using the notation of \refeq{eta.O.x},
$$F((\eta_O)_{t, w} \pp_{\Ax_t} \bigpp x^{\pp t}) = (\eta_P)_{t, w} \pp_{\Ax_t} \bigpp y^{\pp t}$$
by the strong monoidality of $F$.
Consequently, $\Y_P$ is contained in $F(\Y_O)$.
By \csy{\refth{transfer.symmetric.monoidal}\refit{transfer.symmetroidal}}, $\D$ is $F(\Y_O)$-symmetroidal, so we are done.
\xpf

\section{Rectification of algebras over operads}
\label{sect--rectification}

In this section, we use the model structures on modules and algebras over colored operads constructed in the previous section to prove a general operadic rectification result.
Rectification theorems address the following question: given a weak equivalence $P\to Q$ of admissible (symmetric) operads, when are their model categories of algebras Quillen equivalent?

An early rectification for symmetric operads is due to Hinich \cite{Hinich:Homological} in the category $\Ch(\Mod_R)$,
where $R$ is a commutative ring containing $\Q$.
In the same vein, Harper~\cite[Theorem~1.4]{Harper:Monoidal} showed rectification under the assumption that every symmetric sequence is projectively cofibrant.
Lurie \cite[Theorem~4.5.4.7]{Lurie:HA} showed rectification of $\Ei$-algebras to commutative algebras (using the language of $\infty$-operads).
All three results have in common that the model category is required to be freely powered \cite[Definition~4.5.4.2]{Lurie:HA}.

Another class of rectification results applies to symmetric spectra with values in some model category~$\C$.
For individual model categories, such as $\C = \Top$, $\C = \sSet$ and motivic spaces, rectification is due to Elmendorf and Mandell~\cite[Theorem~1.3]{ElmendorfMandell:Rings},
Harper~\cite[Theorem~1.4]{Harper:Symmetric}, and Hornbostel~\cite{Hornbostel:Preorientations}, respectively.
For spectra in an abstract model category $\C$, Gorchinskiy and Guletski\u\i~\cite[Theorem~11]{GorchinskiyGuletskii:Positive} have shown an important special case of symmetric flatness.
We show in \csp{\refsect{symmetricity.2}} that the stable positive model structure on symmetric spectra in (essentially) any model category $\C$ is symmetric flat and give several applications of this fact.

For nonsymmetric operads, Muro \cite[Theorem~1.3]{Muro:Homotopy} has shown a rectification result for a weak equivalence between levelwise cofibrant operads, under similar assumptions to the ones of \refth{rect.colored.operad}.

Our rectification result, \refth{rect.colored.operad}, identifies (symmetric) flatness as a necessary and sufficient condition for the rectification of algebras over (symmetric) colored operads.
It extends the first group of the above-mentioned results since being freely powered is a much stronger condition than being symmetric flat.
It also covers the second group of results since the assumptions of \ref{theo--rect.colored.operad} are satisfied for $\C = \Top$, etc.\ (see \csy{\refchap{examples}}).

We finish this section with \refth{quasicategorical.rectification}, a rectification result relating operadic algebras in the strict sense and in the $\infty$-categorical sense introduced by Lurie.

\theo \label{theo--cow.Coll}
Assume that $\C$ is (symmetric) h-monoidal, symmetric monoidal model category that is (a) strongly admissibly generated, or (b) whose weak equivalences are stable under filtered colimits.
Let $g$ be a weak equivalence in $\ssColl_W$.
\begin{enumerate}[(i)]
\item
\label{item--symmetric.flat.implies.ssColl.flat}
If $g$~is (symmetric) flat in~$\C$ (\refde{summary.symmetricity}),
then $g$~is pseudoflat on the $\ssColl_W$-module $\C^W$, meaning $g \pp b$ is a weak equivalence for any cofibration with cofibrant domain~$b\colon X\to Y$ in~$\C^W$,
where $\pp$ denotes the pushout product of morphisms in $\ssColl_W(\C)$.
\item
\label{item--ssColl.flat.implies.symmetric.flat}
If $g\circ X$ is a weak equivalence for any cofibrant object~$X$ in~$\C^W$, then $g$~is (symmetric) flat in~$\C$,
provided that the coproduct functor reflects weak equivalences and that $\C$ is tractable.
\end{enumerate}
\xtheo

\pf
Recall the multi-index conventions explained in \refsect{symmetricity}.
By definition,
$$(g \pp b)_w = \coprod_{s \in \pi_0(\ssSeq_W^\x)} \underbrace{g_{s, w} \pp_{\Ax_s} \bigotimes_{r \in W} b_r^{\t s^{-1}(r)}}_{=:\lambda_s}.\eqlabel{pf1}$$
(sic, \emph{not} $\bigpp_{r \in W} b_r^{\pp s^{-1}(r)}$).
The coproduct is taken in the category $\Ar (\C)$ of morphisms in $\C$ and runs over all isomorphism classes in $\ssSeq_W^\x$ and $\Ax_s$ is the group of automorphisms of some representative of this isomorphism class.
Recall that $\Ax_s$ is trivial in the nonsymmetric case.
In the symmetric case, an isomorphism class amounts to specifying the number of occurrences of each color $r \in W$, and $\Ax_s$ is as in \refeq{Ax.s}.

We define a multi-index $n$ by $n_r := \sharp s^{-1}(r)$ and set $m_k:=\Sigma_n\cdot_{\Sigma_{n-k}\times\Sigma_k}X^{\otimes n-k}\otimes b^{\pp k}$ for $0\le k\le n$.
By \csy{\refle{combinatorial}}, applied to the composition $\emptyset\To1{}X\To1 b Y$, the map $b^{\t n}$
is the (finite) composition of pushouts of the maps~$m_k$, where $1\le k< n$ and $m_n$ (which is not pushed out).
By \csy{\refpr{pushout.product.pushout.morphisms}, \refle{pp.compositions}}, $\lambda_s$ is therefore the composition of pushouts of
$$g_{s,w} \pp_{\Ax_s} m_k.\eqlabel{pf2}$$

\refit{symmetric.flat.implies.ssColl.flat}~We claim that $\lambda_s$ appearing in \refeq{pf1} is a weak equivalence with h-cofibrant (co)domains.
Recall that an h-cofibrant object~$X$ is such that $\emptyset \r X$ is an h-cofibration.
Weak equivalences with h-cofibrant (co)domains are stable under finite coproducts \cite[Lemma~1.4(a)]{BataninBerger:Homotopy}.
Presenting \refeq{pf1} as the filtered colimit over all finite subsets of the indexing set and using assumption~(b), the claim implies \refit{symmetric.flat.implies.ssColl.flat}.
For assumption (a), we use that the transition maps in the filtered diagram
are cobase changes of morphisms of the form $\emptyset\to\lambda_s$, which in their own turn can be presented as a composition
of maps of the form \refeq{cellular.algebraic}.

To show the claim, we focus on the symmetric case and briefly explain the simpler argument in the nonsymmetric case.
By \csy{\refle{preparation.flat.weakly.sat}} (more precisely, replace $\pp$ by $\pp_{\Ax_s}$ there),
for $\lambda_s$ to be a weak equivalence it is enough to show that the maps in \refeq{pf2} are weak equivalences and that $\text{(co)dom}(g_{s,w}) \t_{\Ax_s} m_k$ is an h-cofibration.
The former holds by symmetric flatness, the latter holds by symmetric h-monoidality, using in both cases the cofibrancy of the (co)domains of~$b_r$.

We now show that $\text{(co)dom}(\lambda_s)$ is an h-cofibrant object.
Writing $g_{s,w}\colon A \r B$, this is clear for $\codom(\lambda_s) = B \t_{\Ax_s}Y^{\otimes n}$, which is h-cofibrant by symmetric h-monoidality,
using the cofibrancy of~$Y_r$.
For the domain of~$\lambda_s$, we first observe that $B \t_{\Ax_s}X^{\otimes n}$ is h-cofibrant.
The map from this object to $\dom(\lambda_s)$ is a cobase change of the map $A \t_{\Ax_s} b^{\otimes n}$.
Again using the above filtration, this map is a composition of pushouts of the maps
$A \t_{\Ax_s} m_k$, which are h-cofibrations by symmetric h-monoidality, using the cofibrancy of $X$.
Since h-cofibrations are stable under pushout and composition \cite[Lemma~1.3]{BataninBerger:Homotopy}, this shows the claim.

\refit{ssColl.flat.implies.symmetric.flat}~First, observe that $g\pp b$ is a weak equivalence for any cofibration with cofibrant source~$b\colon X\to Y$ in~$\C^W$.
Indeed, it suffices to show that $A\circ b$ is an h-cofibration, where $A = \dom(g)$, which follows from symmetric h-monoidality and stability of h-cofibrations under colimits of chains \csy{\refle{i.cofibrations}\refit{i.cofibrations.weakly.sat}}.
Indeed, in this case the pushout of $A \circ b$ along $g \circ X$ is a homotopy pushout since $\C$ is left proper, so that $g \pp b$ is a weak equivalence by the 2-out-of-3 axiom.
The coproduct in~\refeq{pf1} is a weak equivalence, hence so are the $\lambda_s$ because the coproduct functor reflects weak equivalences.
Now we use the filtration~\refeq{pf2} and show by induction on~$n$
that the map $g_{s,w}\pp_{\Ax_m\times\Ax_s}(X^{\otimes m}\otimes b^{\pp n})$ in the definition of symmetric flatness
is a weak equivalence for any cofibration~$b$ with cofibrant source~$X$ and any~$m\ge0$.
The case $m=0$ then gives the symmetric flatness of $g$ relative to $b$.
\ppar
The case $n=0$ is true by assumption (recall that $X$ is assumed to be cofibrant).
For $n\ne0$ consider the filtration~\refeq{pf2} (tensored with $X^{\otimes m}$) of the map $g_{s,w}\pp_{\Ax_m\times\Ax_s}X^{\otimes m}\otimes b^{\otimes n}$,
which is a weak equivalence by assumption (extended to morphisms as explained in the previous paragraph).
For $k\ne n$ the term~$g_{s,w}\pp_{\Ax_m\times\Ax_s}X^{\otimes m}\otimes m_k=g_{s,w}\pp_{\Ax_m\times\Ax_{n-k}\times\Ax_k}X^{\otimes m+(n-k)}\otimes b^{\otimes k}$
is a weak equivalence by the inductive assumption,
and the argument in the previous part shows that its cobase change is a weak equivalence.
Thus, the remaining map in the filtration, $g_{s,w}\pp_{\Ax_m\times\Ax_s}X^{\otimes m}\otimes b^{\pp n}$ (we set $k=n$), is also a weak equivalence, as desired.
\ppar
We have established the symmetric flatness relative to the class of cofibrations with cofibrant source.
Trac\-ta\-bility and the weak saturation property for symmetric flatness \csy{\refth{symmetric.weakly.saturated}\refit{symmetric.flat.weakly.sat}} imply the full symmetric flatness property.
\xpf

\rema \label{rema--sequences.symmetric.flat}
In the situation of \refth{cow.Coll}, similar arguments show that for any weak equivalence~$f$ in $\sColl_W(\C)$ and any cofibrant object $B \in \sColl_W(\C)$, $f \circ B$ is a weak equivalence.
For simplicity of notation, we only consider the uncolored case: then $B = \coprod_{n \ge 0} G_n(A_n)$, where $G_n$ places $A_n$ in degree $n$.
Using the fact that $\circ$ preserves filtered colimits in its second variable and the stability of weak equivalences in~$\C$, hence~$\sColl_W(\C)$,
under filtered colimits, we may assume that $B$ is concentrated in finitely many degrees.

So, let $B=\coprod_{i=1}^k G_{n_i}(A_i)$ (finite coproduct), where $A_i \in \Sigma_{n_i} \C$ is a projectively cofibrant object.
The standard formula for multinomial coefficients takes the following form, where $A_i \in \Sigma_{n_i} \C$, $i = 1, \ldots, k$, $k\ge0$:
$$G_m(f) \circ \left (\coprod_i G_{n_i} (A_i) \right) = \coprod G_{n m} \left (\Sigma_{n m} \cdot_{\Sigma_m \rtimes \Sigma_n^m} f \t A^{\t m} \right).$$
The coproduct runs over all partitions $m = \sum_{i=1}^k m_i$.
The multi-index $(m_1, \ldots, m_k)$ will also be denoted by~$m$ and likewise for~$n$.
In line with the notation in \refsect{symmetricity}, we write $mn = \sum m_i n_i$, and $\Sigma_n^m := \prod \Sigma_{n_i}^{\x m_i}$.
The notation $m n$, $\Sigma_m$ and $\Sigma_n^m$ is understood as in \csy{\refde{multi}}.
Moreover, $A^{\t m}$ stands for $\bigotimes_i A_i^{\t m_i}$.
By \csy{\refle{dirty}}, there is an isomorphism of objects in $\C$ (i.e., disregarding the action of $\Sigma_{nm}$),
$$\Sigma_{n m} \cdot_{\Sigma_m \rtimes \Sigma_n^m} f \t A^{\t m} \cong
\left (f \t_{\Sigma_{m'}} A^{\t m'} \right) \t \left ({\Sigma_{n m''}\over\prod \Sigma_{m''} \rtimes \Sigma_n^{m''}} \cdot A^{\t m''} \right ).$$
Here, $m'$ is the subindex of $m$ consisting of the indices $m_i$ for which $n_i = 0$, whereas $m''$ denotes the subindex of~$m$ consisting of the remaining indices.
As above, $\Sigma_{n m''} := \prod \Sigma_{n_j m''_j}$, etc.
The right factor involving the $A_j$ is cofibrant in $\C$ by the pushout product axiom.
The left factor is a weak equivalence by the symmetric flatness of $\C$.
Our claim now follows from the (nonsymmetric) flatness.
\xrema

The following theorem addresses the question of \emph{Quillen invariance} \cite[Definition~3.11]{SchwedeShipley:Equivalences},
also referred to as {\it rectification}, {\it rigidification}, or \emph{strictification}, i.e., when a weak equivalence of (admissible) operads induces a Quillen equivalence of algebras.

\theo \label{theo--rect.colored.operad}
Suppose that $\C$ is a tractable (symmetric) h-monoidal symmetric monoidal model category such that (a) weak equivalences are stable under filtered colimits
or (b) $\C$ is strongly admissibly generated.
Given a map $f\colon O\to P$ of admissible (symmetric) $W$-colored operads in~$\C$,
the induced Quillen adjunction $$f_*:\Alg_O\rightleftarrows\Alg_P:f^*$$
of the corresponding categories of algebras is a Quillen equivalence if and only if
$f\circ A$ is a weak equivalence for any cofibrant object~$A$ in~$\C^W$.
This condition is satisfied if $f$ is (symmetric) flat in $\C$ (a sufficient condition is given by \refle{projectively.cofibrant.symmetric.flat}).
If the coproduct functor reflects weak equivalences (e.g., the model category is pointed, or we work with simplicial sets or topological spaces), then the opposite is true:
if the above adjunction is a Quillen equivalence, then $f$ (more precisely, its individual levels) is (symmetric) flat in $\C$.
\xtheo

\pf
The adjunction exists by \refth{O.Mod}\refit{pullback.pushforward.O.Mod}.
It is a Quillen adjunction since~$f^*$ preserves (acyclic) fibrations.
\ppar
The flatness condition is necessary because the map~$f \circ A$ is the map~$U(X)\to U(f_*X)$ for the cofibrant object~$X=O \circ A$.
The latter map is the underlying map of the (derived) unit map of~$X$, which is a weak equivalence for any Quillen equivalence.
\ppar
Conversely, by \cite[Definition~8.5.20]{Hirschhorn:Model}, we have to show that a morphism $f_*A \To1 a B$ is a weak equivalence if and only if its adjoint,
i.e., the composition $A \To1\eta f^*f_*A \To3{f^*a}f^*B$, is a weak equivalence
for any cofibrant object~$A$ in~$\Alg_O$ and any fibrant object~$B$ in~$\Alg_P$.
The functor~$f^*$ preserves weak equivalences because both model structures are transferred from~$\C^W$,
thus it remains to prove that $\eta$ is a weak equivalence or, equivalently, that the canonical morphism $U(A)\to U(f_*A)$ is a weak equivalence in~$\C^W$.

We are going to prove a stronger assertion by induction:
for any cofibrant object $Y\in\C^W$,
the map $$U(A\sqcup O\circ Y)\to U(f_*(A\sqcup O\circ Y))$$ is a weak equivalence.
We abbreviate $A_Y:=A\sqcup O\circ Y$, $A'_Y:=A'\sqcup O\circ Y$,
$a_Y:=a\sqcup \id_{O\circ Y}$,
so that the above map is $U(A_Y)\to U(f_*(A_Y))$.

Cofibrant objects in~$\Alg_O$ are retracts of cellular objects and the latter are obtained as codomains of transfinite compositions of cobase changes of generating cofibrations,
starting with the initial $O$-algebra.
For this base case $A=\emptyset$ we get the map $f\circ Y\colon O\circ Y\to P\circ Y$, which is a weak equivalence by assumption on $f\colon O\to P$.

\ppar
To prove the inductive step, we consider
the cocartesian square of $O$-, resp.~$P$-algebras
as in \refeq{pushout.O.Alg}:
$$\xymatrix{O \circ X \ar[d]_{O \circ x} \ar[r] & A \ar[d]^a\cr
O \circ X'\ar[r] & A',\cr} \ \ \ \ \ \ \ \
\xymatrix{P \circ X \ar[d]_{P \circ x} \ar[r] & f_*A \ar[d]^{f_*a}\cr
P \circ X'\ar[r] & f_*A',\cr}$$
where $X\to X'$ is a cofibration between cofibrant (by tractability) objects in~$\C^W$.
The maps $O \circ x$ and $a$ are cofibrations in~$\Alg_O$.
Applying the left Quillen functor~$f_*$ to this square gives a cocartesian square of $P$-algebras as depicted on the right side.
Its vertical maps are again cofibrations and all three objects are cofibrant.

To show that the map $$U(A'_Y)=U(A'\sqcup O\circ Y)\to U(f_*(A'\sqcup O\circ Y))=U(f_*(A'_Y))$$ is a weak equivalence for any cofibrant $Y\in\C^W$,
we use the filtration of \refpr{filtration.Harper},
which shows that for any color $w \in W$, the map $U(a_Y)_w \in \C$ is a transfinite composition of cobase changes $b_s\colon B_s\to B'_s$ of morphisms
$$\alpha_{O,x,A_Y,s} := \Env(O, A_Y)_{s, w} \t_{\Ax_s} \bigpp_{r \in W} x_r^{\pp s^{-1}(r)},\quad s\colon I \r W \in \ssSeq_W^\x,\quad I \ne \emptyset.$$
Likewise, $U(f_*a_Y)_w$ is a transfinite composition of cobase changes $c_s\colon C_s\to C'_s$ of morphisms
$$\alpha_{P,x,f_*A_Y,s} := \Env(P, f_*A_Y)_{s, w} \t_{\Ax_s} \bigpp_{r \in W} x_r^{\pp s^{-1}(r)},\quad s\colon I \r W \in \ssSeq_W^\x,\quad I \ne \emptyset.$$
Since weak equivalences in $\C$ are closed under ordinal-indexed colimits or $\C$ is strongly admissibly generated, it suffices to show that
for any~$s$, if the map $B_s\to C_s$ is a weak equivalence, then so is the induced map $B'_s\to C'_s$.
These maps fit into the following cocartesian squares:
$$
\xymatrix{\dom \alpha_{O,x,A_Y,s}    \ar[d]_{\alpha_{O,x,A_Y,s}}    \ar[r] & B_s \ar[d]^{b_s}\cr \codom \alpha_{O,x,A_Y,s}    \ar[r] & B'_s,\cr}
\qquad
\xymatrix{\dom \alpha_{P,x,f_*A_Y,s} \ar[d]_{\alpha_{P,x,f_*A_Y,s}} \ar[r] & C_s \ar[d]^{c_s}\cr \codom \alpha_{P,x,f_*A_Y,s} \ar[r] & C'_s.\cr}
$$
Both squares are homotopy cocartesian squares because the maps $\alpha_{O,x,A_Y,s}$ and $\alpha_{P,x,f_*A_Y,s}$ are h-cofibrations
by (symmetric) h-monoidality of~$\C$.
We also have a natural transformation between these squares induced by the morphism $O\to P$.
Thus, to show that the map $B'_s\to C'_s$ is a weak equivalence, it suffices to show
that the maps
$$\codom \alpha_{O,x,A_Y,s}\to \codom \alpha_{P,x,f_*A_Y,s},\qquad \dom \alpha_{O,x,A_Y,s}\to \dom \alpha_{P,x,f_*A_Y,s}$$
are weak equivalences.
\ppar
The map $$\codom \alpha_{O,x,A_Y,s}\to \codom \alpha_{P,x,f_*A_Y,s}$$ is a retract of the $w$-component of the map
$\Env(O,A_Y)\circ X'\to \Env(P,f_*A_Y)\circ X'$.
Using the universal property of enveloping operads, cf.~\cite[Proposition~4.7]{Harper:Symmetric} and the definition of $A_Y$, we have isomorphisms, for any $Z \in \C^W$:
$$\Env(O, A_Y) \circ Z = A \sqcup (O \circ (Y \sqcup Z))$$
and similarly with $\Env(P, f_* A_Y)$.
Thus, the map above is isomorphic to
$$A\sqcup(O\circ (Y\sqcup X'))\to f_*(A\sqcup (O\circ (Y\sqcup X'))).$$
which is a weak equivalence by the inductive assumption on~$A$.
\ppar
The map $$\dom \alpha_{O,x,A,s}\to \dom \alpha_{P,x,f_*A,s}$$
admits a filtration defined in \csy{\refle{combinatorial}}, generalized to the multicolored case:
the map $$\Env(O, A_Y)_{s, w} \t_{\Ax_s} \bigotimes_{r \in W} X_r^{\otimes s^{-1}(r)} \to\dom \alpha_{O,x,A,s}$$
is the composition of cobase changes of the maps of the form
$$\Env(O, A_Y)_{s, w} \t_{\Ax_J\cdot\Ax_{s'}} \left( x_q^{\pp J} \otimes \bigotimes_{r \in W} X_r^{\otimes s'^{-1}(r)} \right),
\quad I'\subset I,
\quad I'\ne \emptyset,$$
where $q\in W$ is fixed,
$J=s^{-1}(q)\setminus s'^{-1}(q)$
and for all other $q'\ne q$ we have $s^{-1}(q')\setminus s'^{-1}(q')=\emptyset$.
Thus, it suffices to establish that the induced map on codomains
$$
\Env(O, A_Y)_{s, w} \t_{\Ax_J\cdot\Ax_{s'}} \left( {X'_q}^{\otimes J} \otimes \bigotimes_{r \in W} X_r^{\otimes {s'}^{-1}(r)} \right)
\to
\Env(P, f_*A_Y)_{s, w} \t_{\Ax_J\cdot\Ax_{s'}} \left( {X'_q}^{\otimes J} \otimes \bigotimes_{r \in W} X_r^{\otimes {s'}^{-1}(r)} \right)
$$
is a weak equivalence (the induced map on domains, which themselves are domains of pushout products, is shown to be a weak equivalence in the same manner, by induction).
The latter map is a retract of the map
$$\Env(O,A_Y)\circ(X\sqcup X')\to \Env(P,f_*A_Y)\circ(X\sqcup X'),$$
which, using the universal property of $\Env$ as above, is isomorphic to the map
$$A\sqcup (O\circ (Y\sqcup X\sqcup X'))\to f_*A \sqcup (P\circ (Y\sqcup X\sqcup X'))$$
which in turn is a weak equivalence by the inductive assumption on~$A$.
This finishes the cofibrant induction step.

\ppar
Given a transfinite composition $S=\colim S_i$ in~$\Alg_O$, the map $U(S)\to U(f_*S)$ is a weak equivalence if all maps $U(S_i)\to U(f_*S_i)$ are weak equivalences
because $U$ creates filtered colimits and weak equivalences in~$\C^W$ are stable under filtered colimits by assumption (a).
In case (b), we additionally use that the transition maps $U(S_i) \r U(S_{i+1})$ and similarly with $f_* S_i$ are transfinite compositions
of cobase changes of maps of the form in \refeq{cellular.algebraic}, as witnessed by the filtration \refeq{U.a}.
\xpf

\lemm
\label{lemm--projectively.cofibrant.symmetric.flat}
In a tractable symmetric monoidal model category $\C$, any weak equivalence $f\colon O\to P$ between projectively cofibrant symmetric sequences is symmetric flat.
\xlemm

\pf
We have to show that $g:=f\pp_{\Sigma_n}s^{\pp n}$ is a weak equivalence for any finite family of cofibrations~$s$.
The morphism fits into a commutative triangle with $f\otimes_{\Sigma_n}\codom(s^{\pp n})$ and a cobase change of $f\otimes_{\Sigma_n}\dom(s^{\pp n})$ as the other two sides.
One checks that $\Sigma_n^\proj \C \x \Sigma_n^\inje \C \To 7 {- \t_{\Sigma_n} - } \C$ is a left Quillen bifunctor, where $\proj$ and $\inje$ refers to the projective and injective model structure, respectively.
By tractability, the (co)domain~$X$ of $s^{\pp n}$ is injectively cofibrant,
which implies that $f\otimes_{\Sigma_n}X$ is a weak equivalence.
It remains to observe that the above cobase change is also a homotopy cobase change.
Indeed, $O \t_{\Sigma_n} -$ preserves the injective cofibration~$s^{\pp n}$ with injectively cofibrant source.
The third vertex in the pushout diagram, $O \t_{\Sigma_n} \codom s^{\pp n}$, is also cofibrant, so we have a homotopy pushout.
\xpf

\rema
\refth{rect.colored.operad} is also true for modules (as opposed to algebras) over weakly equivalent operads.
This follows from \refre{sequences.symmetric.flat}.
\xrema

\rema \label{rema--module.context.rectification}
Rectification also holds in a slightly more general context (cf.~\refre{module.context}):
$\C$~is a symmetric monoidal model category, $\C'$ is a tractable model category
whose weak equivalences are stable under filtered colimits and
that is a $\C$-algebra (in the symmetric case, a commutative $\C$-algebra).
Finally, suppose $\C'$ is (symmetric) flat as an algebra (respectively, commutative algebra) over $\C$ (again using an obvious extension of \refde{summary.symmetricity}).
Then any weak equivalence of $W$-colored admissible operads $O \r P$ in $\C$ yields a Quillen equivalence of their algebras in $\C'$.
\xrema

We finish this section by establishing a quasicategorical rectification result, which generalizes~\cite[Theorem~4.5.4.7]{Lurie:HA}
to the case of arbitrary symmetric quasicategorical operads (as opposed to just the commutative operad)
and uses conditions that are significantly weaker than freely poweredness.
The following proposition and theorem, as well as the fact that the former is relevant for the latter, were suggested to the first author by Thomas Nikolaus.
Our proofs are quite similar to that of Lurie in~\cite{Lurie:HA},
the most noticeable difference being the usage of notions of strong admissibility and symmetric flatness.
In particular, strong admissibility allows us to give a rather concise proof
of the preservation of cofibrant objects in the following proposition.

\prop
\label{prop--sifted.homotopy.colimits}
Suppose that $\C$~is a $\V$-enriched cofibrantly generated symmetric monoidal model category and $O$~is a symmetric colored operad in~$\V$ that is admissible in~$\C$.
If the unit map $\eta_O\colon 1[1] \to O$~is a cofibration in $\ssColl_W(\C)$
then the forgetful functor $U\colon\Alg_O(\C)\to\C$ creates (i.e., preserves and reflects) homotopy sifted colimits.
\xprop

\rema
We remind the reader that the notion of a sifted homotopy colimit is stronger than that of a sifted colimit.
For example, the reflexive coequalizer diagram is sifted but not homotopy sifted \cite[Remark~4.5.(e)]{Rosicky:Homotopy}.
This is unlike the filtered case, where both notions coincide for ordinary categories.
\xrema

\pf
The proof is similar to the proof of~\cite[Lemma~4.5.4.12]{Lurie:HA}.
The functor~$U$ creates weak equivalences, so the reflection property is implied by the preservation property.
Denote by~$I$ an arbitrary homotopy sifted small category, such as~$\Delta^\op$.
By Harpaz, Nuiten, and Prasma \cite[Lemma~A.0.2]{HarpazNuitenPrasma}, the functor $\iota\colon I\to I^\sqcup$ is homotopy cofinal,
where $I^\sqcup$ is the universal cocompletion of~$I$ under finite coproducts.
(Proof: observe that for any $L\in I^\sqcup$ given by the formal coproduct of $i_1,\ldots,i_k\in I$,
the comma category $L/\iota$ is equivalent to the comma category $L/d$ for the diagonal functor $d\colon I\to I^k$.
The functor~$d$ is cofinal because $I$ is homotopy sifted, hence the nerve of $L/d$ is weakly contractible, so $\iota$ is also cofinal.)
Following Proposition~A.0.1 there, since $I\to I^\sqcup$ is a fully faithful functor, the category~$I$ is a retraction of~$I^\sqcup$.
Thus, any $I$-diagram is a restriction of an $I^\sqcup$-diagram and so preservation of $I^\sqcup$-indexed homotopy colimits
implies the preservation of $I$-indexed homotopy colimits.
Henceforth, we assume that $I$ has finite coproducts.
\ppar
We now claim that if $I$ admits finite coproducts, then the projective model structure on $\Fun(I,\C)$ is a monoidal model category whenever $\C$ is monoidal.
To verify the pushout product axiom, consider generating (acyclic) cofibrations $\Mor(i,-)\otimes u$ and $\Mor(i',-)\otimes v$,
where $i$ and~$i'$ are objects of~$I$,
the functor $\Mor(i,-)\colon I\to\Set$ is the corepresentable functor of~$i$,
and the morphisms $u$ and~$v$ are (acyclic) cofibrations of~$\C$.
The pushout product of $\Mor(i,-)\otimes u$ and $\Mor(i',-)\otimes v$
is $\Mor(i\sqcup i',-)\otimes(u\pp v)$.
(Here we used the existence of coproducts in~$I$
to show that $\Mor(i\sqcup i',-)\cong\Mor(i,-)\times\Mor(i',-)$.)
By the pushout product axiom for~$\C$, the morphism $u\pp v$ is an (acyclic) cofibration in~$\C$.
Thus, $\Mor(i\sqcup i',-)\otimes(u\pp v)$ is an (acyclic) cofibration in~$\C$.
\ppar
We have a (strictly) commuting diagram
$$\cd{\Fun(I,\Alg_O(\C))&\mapright{\colim}&\Alg_O(\C)\cr
\mapdown{V}&&\mapdown{U}\cr
\Fun(I,\C)&\mapright{\colim}&\C,\cr}$$
where $V$ is also a forgetful functor.
Preservation of homotopy colimits means that the diagram commutes up to a weak equivalence after we derive it.
Both $U$ and~$V$ are automatically derived because they preserve weak equivalences.
We endow $\Fun(I,\Alg_O(\C))$ with the projective model structure (with respect to~$I$) and the transferred model structure on $\Alg_O(\C)$, which exists by assumption.
This model structure is the same as the model structure transferred from the projective model structure on $\Fun(I,\C)$,
if we regard $O$ as an $I$-constant operad in $\Fun(I,\C)$.
Indeed, both model structures are transferred twice: once for the functor category, and the other time for operadic algebras, and it does not matter in which order to transfer.
\ppar
The top $\colim$ (hence, also $U\circ\colim$) can be derived by performing a cofibrant replacement in the source category.
If $V$ preserves cofibrant objects, then it can also be derived in this way, which proves the desired commutativity.
To show that $V$ preserves cofibrant objects,
we observe that $V$ can be rewritten as the forgetful functor $\Alg_O(\Fun(I,\C))\to\Fun(I,\C)$.
It preserves cofibrant objects since $O$~is strongly admissible in $\Fun(I,\C)$ by \refpr{cofibrant.strongly.admissible}.
\xpf

We are now ready to state the conditions under which every quasicategorical algebra
over a quasicategorical operad corresponding to a strict colored symmetric operad can be rectified to a strict algebra over the strict operad.
We state the theorem for the simplicial case, because a detailed write-up of quasicategorical operads is only available in this setting,
however, the proof holds more generally as indicated in the remark below.
This extends results of Lurie~\cite[Theorems~4.1.8.4, 4.5.4.7]{Lurie:HA} for the associative operad and commutative operad,
Haugseng~\cite[Theorem~2.16]{Haugseng:Rectification} for arbitrary nonsymmetric operads,
and Hinich~\cite[Theorem~4.1.1]{Hinich:Rectification} for symmetric operads in the case $\C = \Ch(\Mod_R)$.

\theo
\def\HAlg{\mathop{\bf HAlg}\nolimits}
\label{theo--quasicategorical.rectification}
Suppose that $\C$~is a simplicial symmetric monoidal model category and $O$~is a $\C$-admissible simplicial symmetric colored operad.
Denote by $\CO_\C$ and $\CO_{\Alg_O(\C)}$ the full subcategories spanned by the corresponding classes of cofibrant objects.
The canonical comparison functor $$\N(\CO_{\Alg_O(\C)})[\we_{\Alg_O(\C)}^{-1}]\to\HAlg_{\N^\otimes O}(\N(\CO_\C)[\we_\C^{-1}])$$
is an equivalence of quasicategories if and only if the map $\mcr O \r O$ (the levelwise projective cofibrant
replacement of the underlying symmetric sequence of~$O$) is symmetric flat in $\C$ (\refde{summary.symmetricity}).
Here, $\HAlg$ is used in the sense of Definition~2.1.2.7 (denoted by~$\Alg$ there) in Lurie~\cite{Lurie:HA}
and $\N^\otimes O$ denote the operadic nerve of~$O$, as explained in Definition~2.1.1.23 there.
\xtheo

\rema
If $O$ is nonsymmetric, projective cofibrancy can be replaced by injective cofibrancy (tautologically true for simplicial sets)
because we do not have to mod out symmetric group actions.
Thus, the condition of symmetric flatness can be dropped and every nonsymmetric simplicial colored operad admits quasicategorical rectification.
\xrema

\pf
The symmetric sequence~$\mcr O$ can be constructed by taking the levelwise product of the Barratt--Eccles operad~$\Ei$ and~$O$,
which in fact gives us an operad and not just a symmetric sequence.
The individual levels have a free action of the symmetric group and therefore are projectively cofibrant.
(The levels of $O$ are injectively cofibrant, since any simplicial set is cofibrant.)
They are weakly equivalent to those of~$O$ because simplicial sets are flat and every simplicial set is cofibrant.
\ppar
The morphism $\mcr O\to O$ induces an equivalence of the quasicategories of algebras over $\N^\otimes \mcr O$ and $\N^\otimes O$,
and below we will prove that the comparison functor is an equivalence of quasicategories for~$\mcr O$,
so by the 2-out-of-3 property for equivalences of quasicategories
the main statement is equivalent to $\mcr O\to O$ inducing a Quillen equivalence, which by \refth{rect.colored.operad} is equivalent to symmetric flatness of $\mcr O\to O$.
It remains to show that the comparison map is an equivalence of quasicategories when $O$~is levelwise projectively cofibrant.
\ppar
The rest of the proof coincides with the proof of~\cite[Theorem~4.5.4.7]{Lurie:HA} (modified in the obvious fashion for colored operads instead of the commutative operad),
with the following modifications: for the part~(d) (preservation of homotopy colimits of simplicial diagrams) we use \refpr{sifted.homotopy.colimits},
whereas for part~(e) we have to establish that the free (strict) $O$-algebra on a cofibrant object~$C\in\C^W$
is also the free quasicategorical $O$-algebra in the sense of~\cite[Definition~3.1.3.1]{Lurie:HA}.
Using Proposition~3.1.3.13 there this reduces to proving that the free $O$-algebra $O\circ C=\coprod_{n\ge0}O_n\otimes_{\Sigma_n}C^{\otimes n}$
is also the derived free $O$-algebra.
By assumption $O$~is levelwise projectively cofibrant, so the individual terms in the coproduct are cofibrant in~$\C^W$ and compute the corresponding derived tensor product.
Coproducts of cofibrant objects are also homotopy coproducts, which concludes the proof.
\xpf

\rema
The same proof works (and therefore the theorem holds) for enriched quasicategorical operads
as soon as one has the obvious analog of \cite[Proposition~3.1.3.13]{Lurie:HA}.
We refer the reader to the work of Hongyi Chu and Rune Haugseng
on enriched quasicategorical operads for the case of an arbitrary enriching symmetric monoidal quasicategory.
\xrema

\section{Transport of operads and operadic algebras}
\label{sect--transport}

This section gives an answer to the following important question:
When does a Quillen equivalence $\C \rightleftarrows \D$ of symmetric monoidal model categories induce a Quillen equivalence of (symmetric) operads and their algebras?
The first result in this direction, for monoids and modules over monoids, is due to Schwede and Shipley \cite[Theorem~3.12]{SchwedeShipley:Equivalences}.
This was generalized to nonsymmetric operads and their algebras by Muro \cite[Theorem~1.1,~1.5]{Muro:HomotopyII}, \cite{Muro:Corrections}.
In both statements, the monoidal unit was assumed to be cofibrant.
This assumption, however, is not satisfied in the very interesting stable positive model structure on symmetric spectra \csp{\refth{stable.R.general}}, so we pay special attention to not assuming the cofibrancy of the monoidal unit~$1$.
For example, \refle{cofibrant.replacement}, which governs certain cofibrant replacements, is trivial if $1$ is cofibrant.

\defi
\label{defi--weak.symmetric.monoidal} \cite[Definition~3.6]{SchwedeShipley:Equivalences}
An adjunction between symmetric monoidal categories
$$F:\C \leftrightarrows \D:G \eqlabel{adjunction.F.G2}$$
is a \emph{(symmetric) oplax-lax adjunction} if $G$ is symmetric lax monoidal (see, for example, \cite[Definition~6.4.1]{Borceux:2}).
It is a \emph{weak symmetric monoidal Quillen adjunction} if in addition the oplax structure maps of~$F$ induced from the lax structure of~$G$,
$$\eqalign{F(\mcr 1_\C) & \to 1_{\D},\cr
F(C \t C') & \to F(C) \t F(C')\cr}$$
are weak equivalences for all cofibrant objects $C, C' \in \C$.
\xdefi

In analogy to \refde{well.pointed}, we introduce the following notion.

\defi \label{defi--monoidally.cofibrant}
An object $A$ in a monoidal model category is \emph{well-pointed} if there is a cofibration $1 \r A$. 
\xdefi

As far as their monoidal properties are concerned, well-pointed objects behave like cofibrant objects, as is illustrated by the following lemmas:

\lemm \label{lemm--monoidally.cofibrant.pushout.product}
Let $\C$ be a monoidal model category.
\begin{enumerate}[(i)]
\item \label{item--atB}
If $B$ is well-pointed, then $- \t B\colon \C \r \C$ is a left Quillen functor.
(Thus, well-pointed objects are pseudocofibrant in the sense of Muro~\cite[Appendix~A]{Muro:HomotopyII}.)
\item \label{item--appb}
If $a\colon A \r A'$ and $b\colon B \r B'$ are two cofibrations with well-pointed source, then so is $a \pp b$.
If either $A$ or $B$ is cofibrant, then $a \ppdom b$ is also cofibrant.
\end{enumerate}
\xlemm

\pf
\refit{atB}~Pick a cofibration $\eta\colon 1 \r B$.
For any (acyclic) cofibration $a$, the map $a \t B$ is the composition of a pushout of $a = a \t 1$ and $a \pp \eta$.
Both are (acyclic) cofibrations.

\refit{appb}~By \refit{atB}, $A \t B$ is well-pointed and $a \t B$ and $A \t b$ are cofibrations.
Hence, $a \ppdom b := \dom(a \pp b)$ is well-pointed as well.
If, say, $A$ is cofibrant, then $\emptyset \r A \To4{A\otimes\eta} A \t B \r a \ppdom b$ is a composition of cofibrations.
\xpf

\lemm \label{lemm--cofibrant.replacement}
Let $A$ and $B$ be two cofibrant or well-pointed objects in a tractable monoidal model category
satisfying the unit axiom,
i.e., $\mcr (1) \t C \sim C$ for all cofibrant objects $C$.
Also assume that (a) weak equivalences are stable under filtered colimits or (b) $\C$ is strongly admissibly generated.
Then the following map is a weak equivalence:
$$\mcr (A) \t \mcr (B) \r A \t B.$$
\xlemm

\pf
If $A$ and $B$ are cofibrant, the claim is clear.
We now show the statement if $B$ is cofibrant and $A$ is well-pointed.

The cofibration $1 \r A$ is a retract of a transfinite composition of maps $A_0 = 1 \r \cdots \r A_\infty = A$,
where each $a_n\colon A_n \r A_{n+1}$ is the pushout of a generating cofibration $s\colon S \r S'$.
We write $E_n\colon s \r a_n$ for the pushout square.
The functor $- \t B$ is a left Quillen functor by \refle{monoidally.cofibrant.pushout.product}\refit{atB}.
In particular, it preserves cofibrations, so that $E_n \t B$ is a pushout of a cofibration between cofibrant objects along a map with cofibrant target $A_n \t B$ (which holds by induction, starting with $A_0 \t B = B$).
Hence, it is a homotopy pushout square.
Similarly, $\mcr(E_n)$ is a pushout one of whose legs is a cofibration, and all objects in the square are cofibrant.
Hence, $\mcr (E_n) \t \mcr (B)$ is also a homotopy pushout square.
In the natural transformation of homotopy pushout squares
$$\mcr (E_n) \t \mcr (B) \To1{} E_n \t B,$$
the two left maps in the depth direction are
$$\mcr (S) \t \mcr (B) \To1\sim S \t B,\eqlabel{Q.S.B}$$
since $\mcr (S) \r S$ is a weak equivalence between cofibrant objects and similarly for~$B$.
(This is the only point where we are using the cofibrancy of~$B$.)
The same works for $S'$.
The third map is
$$\mcr (A_n) \t \mcr (B) \r A_n \t B,\eqlabel{Q.A.B}$$
which by induction on $n$ is a weak equivalence, starting for $n=0$ with the weak equivalence
$$\mcr (1) \t \mcr (B) \sim 1 \t \mcr (B) = \mcr (B) \sim B$$
given by the unit axiom.
Thus, the fourth map in the cube, $\mcr (A_{n+1}) \t \mcr (B) \r A_{n+1} \t B$, is a weak equivalence.
Thus, for all $n < \infty$, \refeq{Q.A.B} is a weak equivalence.
In other words, $\mcr (A_n) \t \mcr (B)$ is a cofibrant replacement of $A_n \t B$.
Then $\mcr (A_\infty) \t \mcr (B) \sim \colim \mcr (A_n) \t \mcr (B) \sim \colim A_n \t B = A_\infty \t B$,
using that weak equivalences are stable under filtered colimits by assumption and the preservation of filtered colimits by $\t$.
In case~(b), we additionally use that the transition maps are cobase changes of generating cofibrations tensored with a fixed object, hence in the class \refeq{cellular.algebraic}.
We have shown the claim if $B$ is cofibrant.

If $B$ is merely well-pointed, we run the same argument again, noting that for a cofibrant object~$S$, the weak equivalence
$\mcr (S) \t \mcr (B) \sim S \t B$ used in \refeq{Q.S.B} is a weak equivalence by the previous step.
\xpf

The following variant can be proved using the same technique as \refle{cofibrant.replacement}.
The left properness is used to ensure that the pushouts appearing in the cellular induction are homotopy pushouts.
The details are left to the reader.

\lemm \label{lemm--monoidally.cofibrant.flat}
Let $A$ be a cofibrant or well-pointed object in a flat left proper tractable monoidal model category $\C$
whose weak equivalences are stable under filtered colimits.
Then $A \t -$ preserves weak equivalences.
\xlemm

The following lemma of Berger and Moerdijk may be called an \emph{equivariant pushout product axiom}.

\lemm \cite[Lemma~2.5.3]{BergerMoerdijk:Boardman}.
\label{lemm--group.pushout.product}
Let $1 \r \Gamma_1 \r \Gamma \r \Gamma_2 \r 1$ be a short exact sequence of finite groups.
Then, for a monoidal model category $\C$,
$$\t\colon \Gamma_2^\proj \C \x \Gamma^{\proj'} \C \r \Gamma^\proj \C$$
is a left Quillen bifunctor.
Here, $\Gamma^{\proj'} \C$ denotes the model structure on $\Gamma \C$ whose cofibrations are $\Gamma_1$-projective cofibrations.
\xlemm

\theo \label{theo--transport}
Suppose $F : \C \leftrightarrows \D: G$ is a weak symmetric monoidal Quillen adjunction (\refde{weak.symmetric.monoidal})
between tractable symmetric monoidal model categories such that (a) weak equivalences are stable under filtered colimits or (b) $\C$ is strongly admissibly generated.
Also suppose that both $\C$ and $\D$ are either left proper or their monoidal unit is cofibrant.
\begin{enumerate}[(i)]
\item \label{item--adjunction.Oper}
Suppose that the transferred model structures on the categories $\ssOper_W(\C)$ and $\ssOper_W(\D)$ exist.
(See \refco{operads} for a sufficient condition.)
Then there is a Quillen adjunction of the categories of (symmetric) operads
$$F^\ssOper: \ssOper_W (\C) \rightleftarrows \ssOper_W (\D) : G. \eqlabel{adjunction.Oper}$$
It is a Quillen equivalence if $(F, G)$ is a Quillen equivalence.
\item \label{item--adjunction.Alg}
For any admissible (symmetric) operad $O$ in $\C$, there is a Quillen adjunction
$$F^\Alg : \Alg_O (\C) \rightleftarrows \Alg_{F^\ssOper(O)} (\D): G.\eqlabel{adjunction.Alg}$$
It is a Quillen equivalence if $(F, G)$ is a Quillen equivalence and $O$ is a cofibrant operad.
\item \label{item--adjunction.Alg2}
If $P$ is an admissible (symmetric) operad in~$\D$ such that $G(P)$ is also admissible, there is a Quillen adjunction
$$F_\Alg: \Alg_{G(P)}^\C\rightleftarrows \Alg_P^\D: G.\eqlabel{adjunction.Alg.P}$$
It is a Quillen equivalence if $(F, G)$ is a Quillen equivalence, $P$ is fibrant, and $\C$ and $\D$ admit rectification of (symmetric) operads.
\end{enumerate}
\xtheo

\pf
Since $G$ is symmetric lax monoidal, it induces a lax monoidal adjunction
$$F: (\ssColl_W \D, \circ) \r (\ssColl_W \C, \circ): G.\eqlabel{adjunction.nsColl.F.G}$$
In particular, $G$ preserves monoids, i.e., (symmetric) operads.
This defines the right adjoint in \refeq{adjunction.Oper}.
The right adjoint in \refeq{adjunction.Alg} sends an $F^\ssOper(O)$-algebra $B$ to $G(B)$, which is an $O$-algebra via
$$O \circ G(B) \r GF^\ssOper(O) \circ G(B) \r G(F^\ssOper(O) \circ B) \r G(B).$$
The left adjoints exist by \cite[Theorem~4.5.6]{Borceux:2}.
Moreover, the right adjoints are Quillen right adjoints since (acyclic) fibrations are again created by the forgetful functors.

We now establish the advertised Quillen equivalences.

\refit{adjunction.Oper}~We have to show that for any cofibrant operad $O$, the natural map
$$\phi_O\colon F(\mcr (U(O))) \r U(F^\ssOper(O))$$
is a weak equivalence.
In this case, we have the following chain of equivalent statements for any cofibrant operad $O \in \ssOper_W(\C)$ and any fibrant operad $P \in \ssOper_W(\D)$, which implies the Quillen equivalence \refeq{adjunction.Oper}:
$$\eqalign{
F^\ssOper(O) \sim P & \Leftrightarrow UF^\ssOper(O) \sim U(P) \cr
& \Leftrightarrow F(\mcr (U(O)) \sim U(P) \cr
& \Leftrightarrow \mcr (U(O)) \sim G(U(P)) = U(G(P)) \cr
& \Leftrightarrow U(O) \sim U(G(P)) \cr
& \Leftrightarrow O \sim G(P).\cr}$$

The cellular induction starts with the initial operad $O = 1_\C[1]$, for which $F^\ssOper(O) = 1_\D[1]$.
Thus, $\phi_{1[1]}$ is a weak equivalence by the weak monoidality of $F$.

Using the notation of \refpr{filtration.operads}, we now consider a pushout of operads along a map $\Free(x)$, where $x$ is a cofibration in $\sColl_W(\C)$.
We will show that $\phi_{O'}$ is a weak equivalence provided that $\phi_O$ is one.

Applying $F\mcr $ to the filtration (see \refpr{filtration.operads})
$$U (o) \colon O^{(0)} := U(O) \r \cdots \r O^{(\infty)} := U(O')$$
gives the front face of the following commutative cube in $\Ax_s \D$.
The back face is part of the filtration
$$U (\tilde o) \colon \tilde O^{(0)} := UF^\ssOper(O) \r \cdots \r \tilde O^{(\infty)} := UF^\ssOper(O')$$
associated to the pushout of operads in $\D$ that is obtained by applying the left adjoint $F^\ssOper$ to \refeq{pushout.operads}:
$$\xymatrix{
\Free(\tilde X) := F^\ssOper(\Free(X)) \ar[d]_{\Free(\tilde x)} \ar[r] & \tilde O := F^\ssOper (O) \ar[d]^{\tilde o} \\
\Free(\tilde X') := F^\ssOper(\Free(X')) \ar[r] & \tilde O' := F^\ssOper(O').
}$$
Here and below, the notation $\tilde ?$ indicates the object or morphism that is obtained by considering the data in the filtration of $\tilde o := F^\ssOper(o)$.
For example, $\tilde X := F(X)$ and similarly for $X'$, $x$.
The coproduct runs over all isomorphism classes of marked trees $T$ in $\ssTree_{s,w}^{(k+1)}$.
$$\xymatrix{
&\!\!\!\!\!\coprod_{T} \Ax_s \cdot_{\Aut T} \tilde x^*(T) \ar'[d][dd] \ar[rr]& &\tilde O^{(k)}_{s,w} \ar[dd] \\
F\mcr \left (\coprod_{T} \Ax_s \cdot_{\Aut T} x^*(T) \right ) \ar[ur]^{*} \ar[dd] \ar[rr] & &F\mcr (O^{(k)}_{s,w}) \ar[dd] \ar[ur]^{r^{(k)}} \\
&\!\!\!\!\!\coprod_{T} \Ax_s \cdot_{\Aut T} \tilde x(T) \ar'[r][rr] & & \tilde O^{(k+1)}_{s,w} \\
F\mcr \left (\coprod_{T} \Ax_s \cdot_{\Aut T} x(T) \right) \ar[rr] \ar[ur]^{**} & &F\mcr (O^{(k+1)}_{s,w}) \ar[ur]^{r^{(k+1)}}\\}$$
At this point (and only here), we use the assumption that $\D$ is either left proper or its monoidal unit is cofibrant:
in the former case, any pushout along a cofibration is a homotopy pushout.
In the latter case, $\tilde O_{s,w} = \tilde O^{(0)}_{s,w}$ is cofibrant for all $(s,w)$ by \refle{cofibrancy.preserved}\refit{special.case.3} and therefore by induction the same is true for $\tilde O^{(k)}_{s,w}$.
Hence, the pushout above is again a homotopy pushout.
Likewise, the front square is a homotopy pushout, since $F\mcr (-)$ preserves those.
Thus, $r^{(k+1)}$ is a weak equivalence if $r^{(k)}$, $*$, and~$**$ are ones.
The map $r^{(k)}$ is a weak equivalence by induction on $k$, starting with
$$r^{(0)} \colon F\mcr (O^{(0)}_{s,w})=F\mcr (U(O)_{s,w}) \r \tilde O^{(0)}_{s,w} = UF^\ssOper(O)_{s,w},$$
which is the $(s,w)$-level of $\phi_O$, which is a weak equivalence by the cellular induction on $O$.
It remains to show that the maps $*$ and $**$ are weak equivalences.

Let $T \in \ssTree_{s,w}$ be any tree.
By induction on the height of $T$, we prove the following claims:
\begin{enumerate}[(A)]
\item
\label{item--claim.A}
The map $\epsilon(T)$ is a cofibration in $(\Aut T)^\proj \C$ with cofibrant or well-pointed domain (\refde{monoidally.cofibrant}).
The domain is cofibrant for all trees except (possibly) for the tree $T_w^- := (\mathord{\mathop\to\limits^w}\mathord{\mathop\bullet\limits^-}\mathord{\mathop\to\limits^w}) \in \ssTree_{w,w}^{(0)}$, which consists of a single nonmarked vertex with input edge and root edge colored by $w$.
In particular, $\epsilon(T)$ is a cofibration with cofibrant domain for all $T \in \ssTree_{s,w}^{(k+1)}$ with $k \ge 0$.
(These are the trees appearing in the cubical diagram above.
In order to perform the induction, we also need to consider $T \in \ssTree^{(0)}_{s,w}$.)
\item
\label{item--claim.B}
There are weak equivalences in $\Ar (\C)$ (i.e., both source and target of the morphisms are weakly equivalent)
$$F\mcr (\epsilon(T)) \r \tilde \epsilon(T).$$
\end{enumerate}

Let $(t,w) := \val (r(T))$ be the valency of the root $r(T)$ of $T$.
If $T$ consists of a single vertex $r(T)$ (with an outgoing root edge and finitely many input edges), then $t=s$ and
$$\epsilon(T) = \epsilon(r(T)) =\cases{
(\eta_O)_{(t,w)},&if the root~$r(T)$ is not marked; \cr
x_{(t,w)},&if the root~$r(T)$ is marked.}$$
Both are cofibrations in $\Ax_t(\C) (= \Aut(T) \C)$, the former by \refle{cofibrancy.preserved}\refit{forget.Oper}.
Since $X = \dom(x)$ is cofibrant by quasitractability, the source of $\epsilon(T)$ is well-pointed for $(T =) r(T) = T_w^-$ and cofibrant else.
This shows claim \refit{claim.A}.

For claim \refit{claim.B}, we note that $F\mcr (U(\eta_O))$ is weakly equivalent to $\eta_{\tilde O}$ by the unit part of the weak monoidality of $F$ and the cellular induction on $O$.
To show $F\mcr (u) \sim \tilde u$, we consider the pushout square in $\ssColl_W(\C)$, denoted by~$E$:
$$\xymatrix{
X \ar[r] \ar[d]^x & U(O) \ar[d]^u\cr
X' \ar[r] & U(O) \sqcup_X X'.\cr
}$$
It is a homotopy pushout square in all degrees: for unit degrees, the left vertical map is $\id_\emptyset$ and for nonunit degrees $O_{s,w}$ is ($\Ax_s$-projectively) cofibrant (and $x_{s,w}$ is a cofibration).
Applying $F\mcr $ to $E$ gives a homotopy pushout square in $\ssColl_W(\D)$.
The square $\tilde E$ in $\sColl_W(\D)$ obtained by replacing $X$, $X'$, and $O$ by their $\tilde ?$-counterparts is also a homotopy pushout square.
By cellular induction $F\mcr U(O) \sim U \tilde O$.
Of course $F\mcr (X) \sim \tilde X (= F(X))$ by the cofibrancy of~$X$ (using the quasitractability of~$\C$) and similarly for~$X'$.
We obtain the desired weak equivalence
$$F(\mcr (U(O) \sqcup_X X')) \sim U(\tilde O) \sqcup_{\tilde X} \tilde X'$$
and hence claim \refit{claim.B} for the tree $T$ consisting of a single (marked or unmarked) vertex.

We now perform the induction step.
We may assume that $T$ has at least two vertices.
By definition,
$$\epsilon(T) = \epsilon(r(T)) \pp \underbrace{ \bigpp_i \epsilon(T_i)^{\pp t_i} }_{=: \epsilon'(T)}.$$

Recall that a map $f$ in a model category $\C$ is a cofibration with cofibrant source if and only if it is a cofibrant object in $\Ar(\C)$, i.e., $\id_\emptyset \r f$ is a cofibration.
Likewise, $f$ is a cofibration with well-pointed source if and only if there is a cofibration $\id_1 \r f$ in $\Ar(\C)$.

We write $\epsilon(r(T))\colon V \r W$ and $\epsilon'(T)\colon e^*(T) \r e(T)$.
Let $\val (r(T)) = (s, w)$.
As was noted above, $\epsilon(r(T))$ is a cofibration in $\sColl_W(\C)$.
Its domain $V_{s, w} := \dom(\epsilon(r(T))_{s,w})$ is well-pointed in $\Ax_s \C$ if $T$ is of the form
$(T_1 \mathord{\mathop\to\limits^w}\mathord{\mathop\bullet\limits^-}\mathord{\mathop\to\limits^w})$, where $T_1$ is the subtree of the root vertex.
In this case, we abusively write $r(T) = T_w$.
In all other cases, $V_{s,w}$ is cofibrant.
Hence, $\id_1 \r \epsilon(r(T))$ (respectively, $\id_\emptyset \r \epsilon(r(T))$) is a cofibration in $\Ar(\Ax_s \C) = \Ax_s \Ar(\C)$.
By induction on $T$, $\epsilon(T_i)$ is an $\Aut (T_i)$-projective cofibration whose source is well-pointed (if $T_i = T_w$) and cofibrant (otherwise).
Again, we reinterpret this in terms of cofibrations in $\Ar(\Aut(T_i) \C)$.

We now consider four cases:
\begin{enumerate}
\item
$r(T) \ne T_w^-$, at least one $T_i \ne T_w^-$:
By \refle{group.pushout.product}, applied to $\Ar(\C)$ (with the pushout product), the map
$$(\id_\emptyset \r \epsilon(r(T))) \pp (\id_\emptyset \r \epsilon'(T)) = (\id_\emptyset \r \epsilon(r(T)) \pp \epsilon'(T)) = (\id_\emptyset \r \epsilon(T))$$
is a cofibration in $\Ar(\Aut(T) \C)$ in this case, i.e., $\epsilon(T)$ is a cofibration with cofibrant source.
\item
$r(T) \ne T_w^-$, all $T_i = T_w^-$:
Then
$$(\id_\emptyset \r \epsilon(r(T))) \pp (\id_1 \r \epsilon'(T)) = (\id_\emptyset \r \epsilon(r(T)) \pp \epsilon'(T)) = (\id_\emptyset \r \epsilon(T))$$
is a cofibration in $\Ar(\Aut(T) \C)$.
\item
Similarly for $r(T) = T_w^-$, $T_1 \ne T_w^-$.
\item \label{item--impossible.case}
$r(T) = T_w^-$, $T_1 = T_w^-$: By definition of the trees in $\ssTree_{s,w}$, any internal edge contains at least one marked vertex.
Thus, this tree does not lie in $\ssTree_{s,w}$ unless $T_1$ is empty, in which case we have shown the claim above.
\end{enumerate}
This shows claim \refit{claim.A}.

We now show \refit{claim.B}.
We may assume that $T$ consists of at least two vertices.
Consider the diagram $E$ whose left square is by definition cocartesian,
$$\xymatrix{
V_{t,w} \t e^*(T) \ar[rr]^{V_{t,w} \t \epsilon'(T)} \ar[d]_{\epsilon(r(T))_{t,w} \t e^*(T)} & &
V_{t,w} \t e(T) \ar[d] \ar[drr] \\
W_{t,w} \t e^*(T) \ar[rr] & &
P \ar[rr]_{\epsilon(r(T))_{t,w} \pp \epsilon'(T)} & &
W_{t,w} \t e(T).
}
\eqlabel{5.terms.operads}$$
We claim that the left pushout square is a homotopy pushout.
By \refle{monoidally.cofibrant.pushout.product}\refit{atB}, both the left vertical and the top horizontal maps are cofibrations (in $\C$, say),
hence the claim is clear if $V_{t,w} \t e^*(T)$ is cofibrant, because in this case the above pushout diagram is cofibrant as a diagram.
By the above, $V_{t,w}$ and $e^*(T)$ are either cofibrant or well-pointed.
Again using \refle{monoidally.cofibrant.pushout.product}, the only way that $V_{t,w} \t e^*(T)$ is only well-pointed is that both $V_{t,w}$ and $e^*(T)$ are well-pointed.
By the above, the first only happens for $r(T) = T_w^-$ and the second happens only if all $T_i = T_w^-$.
As was noted in Case \refit{impossible.case}, this means $T=(\mathord{\mathop\to\limits^w}\mathord{\mathop\bullet\limits^-}\mathord{\mathop\to\limits^w}\mathord{\mathop\bullet\limits^-}\mathord{\mathop\to\limits^w})$, which is excluded.

We have weak equivalences
$$\eqalign{F\mcr (V_{t,w} \t e^*(T)) & \sim F(\mcr V_{t,w} \t \mcr e^*(T)) \cr
& \sim F\mcr (V_{t,w}) \t F\mcr (e^*(T)) \cr
& \sim \mcr (\tilde V_{t,w}) \t \mcr (\tilde e^*(T)) \cr
& \sim \tilde V_{t,w} \t \tilde e^*(T).\cr}$$
The first equivalence holds by \refle{cofibrant.replacement}, which gives a weak equivalence between cofibrant objects
$$\mcr (V_{t,w} \t e^*(T)) \sim \mcr (V_{t,w}) \t \mcr (e^*(T))$$
since both $V_{t,w}$ and $e^*(T)$ are cofibrant or well-pointed.
The second equivalence holds by weak monoidality of $F$.
The third equivalence follows from Brown's lemma and the equivalences $F\mcr (V_{t,w}) \sim \tilde V_{t,w}$ and $F\mcr (e^*(T)) \sim \tilde e^*(T)$.
The last weak equivalence holds by \refle{cofibrant.replacement}, again using the (monoidal) cofibrancy of $\tilde V_{t,w}$ and $\tilde e^*(T)$.
The same is also true for $W_{t,w}$ and/or $e(T)$ instead.

We now apply $F\mcr $ to the diagram $E$ in \refeq{5.terms.operads}.
On the other hand, we consider the diagram $\tilde E$ obtained by replacing $V_{t,w}$ by $\tilde V_{t,w}$, etc.
There is a map of diagrams $F\mcr (E) \r \tilde E$.
By the above, all individual maps in this morphisms of diagrams are weak equivalences, except (a priori) for
$$F\mcr \left( P \right ) \r
\tilde P.$$
However, since the left squares of $F\mcr (E)$ and $\tilde E$ are homotopy pushout squares, this remaining map is also a weak equivalence.
Therefore, $F\mcr (E) \sim \tilde E$.
In particular, we get the requested weak equivalence in $\Ar(\C)$
$$F\mcr (\epsilon(T)) \sim \tilde \epsilon (T).$$
This finishes the induction step (with respect to the tree $T$).
We have shown that the individual summands in the maps $*$ and $**$ are weak equivalences.

The coproducts appearing in the left face of the cube above are homotopy coproducts,
since for all $T \in \ssTree_{s,w}^{(k+1)}$ ($k \ge 0$), the terms $\Ax_t \cdot_{\Aut T} x^*(T)$ and similarly for $x(T)$ are $\Ax_t$-projectively cofibrant by Claim \refit{claim.A}.
This implies that the maps $*$ and $**$ themselves are weak equivalences and therefore finishes the induction step with respect to the cellular induction by $O$.

For a cellular filtration of~$O_\infty$ by operads~$O_i$ such that $\phi_{O_i}$ is a weak equivalence for all $i < \infty$,
the same is true for $i = \infty$ using that $U$ preserves filtered colimits and assumption~(a).
In case (b), we also use that the transition maps $\text{(co)dom}(\phi_{O_i}) \r \text{(co)dom}(\phi_{O_{i+1}})$ lie in \refeq{cellular.algebraic}, by \refeq{filtration.operad}.

\refit{adjunction.Alg}~For any cofibrant $O$-algebra $A$,
we have the following chain of canonical isomorphisms and weak equivalences, which as above shows the requested Quillen equivalence:
$$\eqalignno{
U(F^\Alg(A)) & = \Env(F^\ssOper (O), F^\Alg(A))_0 &\eqlabel{proof.adjunction.Alg} \cr
&= F^\ssOper(\Env(O, A))_0 \cr
&\sim F(\mcr (\Env(O, A))_0) \cr
&\sim F(\Env(O, A)_0) \cr
&= F(U(A)) \cr
&\sim F\mcr (U(A)).\cr}$$
The last (and similarly the first) canonical isomorphism is \refpr{enveloping.operad}\refit{O.A.initial}.
The second isomorphism comes from a natural isomorphism of functors
$$\Env(F^\ssOper (-), F^\Alg(*)) = F^\ssOper(\Env(-, *))$$
since both expressions are the left adjoint to $\ssOper_W(\D) \r \Pairs(\sColl_W(\C))$, $P \mapsto (G(P), G(P)_0)$.
The first weak equivalence was shown in Part \refit{adjunction.Oper}, which is applicable since $\Env(O, A)$ is a cofibrant operad by \refle{cofibrancy.preserved}\refit{Env.preserves.cofibrations}.
The second weak equivalence is given by \refle{cofibrancy.preserved}\refit{forget.Oper}.
The last weak equivalence follows from \refpr{cofibrant.strongly.admissible}.

\refit{adjunction.Alg2}~Let $O \in \ssOper_W(\C)$ be a cofibrant replacement of $G(P)$.
Equivalently, by Part \refit{adjunction.Oper},
$P \sim F^\ssOper (O)$.
By rectification of operads for $\D$, \refit{adjunction.Alg}, and rectification of operads for $\C$, we have the following chain of Quillen equivalences
$$\Alg_P^{\D} \sim \Alg_{F^\ssOper(O)}^{\D} \sim \Alg_O^\C \sim \Alg_{G(P)}^\C.$$
\xpf

\rema
The condition in \refth{transport} that $\C$ and $\D$ have the property that they are either left proper or their monoidal unit is cofibrant is only used to show that pushouts of certain cofibrations with cofibrant domain are homotopy pushouts.
Since being a homotopy pushout only depends on the class of weak equivalences, this also holds, for example, if $\C$ has another model structure with more cofibrations, and the same weak equivalences.
\xrema

If the left adjoint $F$ is in addition symmetric monoidal, we can relax the condition on $O$ in \refth{transport}\refit{adjunction.Alg}.

\coro
\label{coro--transport.F.monoidal}
In the situation of \refth{transport}, suppose in addition that the left adjoint $F$ is strong symmetric oplax monoidal (i.e., the symmetric oplax structural maps $F(C \t C') \r F(C) \t F(C')$ are isomorphisms, so that $F$ is also symmetric lax monoidal).
Let $O$ be any (symmetric) operad in $\C$ such that $U(\eta_O)$ is a cofibration in $\ssColl_W(\C)$.

Then there is a Quillen adjunction
$$F: \Alg_O (\C) \rightleftarrows \Alg_{F(O)} (\D): G,$$
which is a Quillen equivalence if $(F, G)$ is a Quillen equivalence.
\xcoro

\pf
Since $F$ is symmetric monoidal, $U \circ F^\Alg = F \circ U$ (see, for example, \cite[Proposition~3.91]{AguiarMahajan:Monoidal}).
Therefore, only the last weak equivalence in \refeq{proof.adjunction.Alg} requires proof.
By \refpr{cofibrant.strongly.admissible}, the operad $O$ is strongly admissible, i.e., $U(A)$ is cofibrant in $\C$, so that $F(U(A)) \sim F(Q(U(A))$ by Brown's lemma.
\xpf

\numberwithin{subsection}{section}
\numberwithin{equation}{subsection}

\section{Applications}
\label{chap--applications}

This last section contains a few applications to the homotopy theory of enriched categories, ordinary categories, operads, and (monoidal) diagrams.
The strategy is similar for all these applications: enriched categories, say, are algebras over a certain nonsymmetric operad.
Therefore, the admissibility and rectification results of \S\ref{sect--operad.admissible}--\ref{sect--rectification} can be applied.

The list presented here is by no means exhaustive,
other potential applications include monads in model categories, internal categories (and higher internal categories), (higher) spans, etc.
Symmetric operads in symmetric spectra and some applications are studied in~\cite{PavlovScholbach:Spectra}.

In \refchap{applications}, let $\V$ be a symmetric monoidal model category and $\C$ be a $\V$-enriched model category whose weak equivalences are stable under filtered colimits.
Moreover, assume that $\C$ is tractable and either combinatorial or $\V$-admissibly generated.

\subsection{Rectification of $\Ai$- and $\Ei$-monoids}

In this section, we discuss rectification of homotopy coherent versions of monoids and commutative monoids.
We start by giving explicit constructions of two important operads, $\Ai$~and~$\Ei$.

The Barratt--Eccles operad~$\Ei$ can be constructed by taking the associative symmetric operad in sets, applying the functor~$\E$ to it
($\E$ sends a set to a groupoid with the same set of objects and a single morphism between any pair of objects), obtaining a symmetric operad in groupoids,
and then applying the nerve functor, which gives a simplicial operad.
See the paragraph after Corollary~3.5 in Elmendorf and Mandell~\cite{ElmendorfMandell:Rings}.

An identical construction (apply~$\E$ and then take the nerve) produces a model for the operad~$\Ai$,
but the original operad in sets is now the free operad on a single binary operation and a single nullary operation,
so that $\cO_n$ consists of planar rooted trees with $n$~leaves (see, for example,~\cite[\S5.8]{BergerMoerdijk:Axiomatic}).
Alternatively, one can take the free operad generated by a single operation in each arity (which corresponds to the so-called {\it unbiased\/} monoids).

In what follows, we actually do not need to apply the nerve functor, because an operad in groupoids is sufficient for our purposes.
We also note that any category enriched in simplicial sets is automatically enriched in groupoids by applying the nerve functor.
The following propositions are mere specializations of the general theorems on admissibility and rectifiability.
We give explicit statements here due to the importance of these examples.

\prop
If $\C$ is a symmetric h-monoidal and groupoid-enriched
then the category of $\Ei$-al\-ge\-bras in~$\C$ admits a transferred model structure.
Furthermore, if the morphism $\Ei\to\Comm$ is symmetric flat in $\C$ (in particular, if $\C$ is symmetric flat),
then the Quillen adjunction between commutative monoids and $\Ei$-monoids is a Quillen equivalence.

A similar statement for $\Ai$ and $\Ass$ holds if $\C$ is merely h-monoidal and flat.
\xprop

\subsection{Model structures on enriched categories}
\label{sect--enriched.categories}

For a small set $W$, Berger and Moerdijk \cite[1.5.4]{BergerMoerdijk:Resolution} have introduced a nonsymmetric $W\times W$-colored operad in~$\V$ given by
$$\CatOp^\Ass_W(((v_1, v_1'), \ldots, (v_n, v_n')), (v'_0, v'_{n+1}))
=\cases{1_\V,&$v'_i = v_{i+1}$ for all $0\le i\le n$;\cr
\emptyset,&otherwise.\cr}$$
This defines a nonsymmetric operad in $\V$.
Its algebras in $\C^{W \x W}$ are precisely \emph{$\C$-enriched categories with $W$ as objects}.
More generally, given a nonsymmetric operad $O$ in $\V$, one can also consider the nonsymmetric operad $\CatOp^O_W$, which is given by replacing $1_\V=\Ass_n$ in the previous formula by~$O_n$.
Algebras over this operad can be called \emph{$\V$-enriched $O$-twisted categories}.
Typically, $O$ is taken to be $\Ai$.
In this case we speak of \emph{$\V$-enriched $\Ai$-categories}, i.e., composition is not strictly associative, but rather associative up to coherent higher homotopies.

The following lemma is an immediate application of the results on admissibility and rectification.
Up to a minor expository difference (see \refre{explain.assumptions}), the admissibility statement is the same as Muro's \cite[Corollaries 10.4,~10.5]{Muro:Homotopy}.
The rectification result there uses in addition the left properness of $\C$.

\coro
\label{coro--Cat.C.fixed.objects}
If $\C$ is h-monoidal, then all (nonsymmetric) operads in~$\V$ are admissible.
In particular, the operad $\CatOp^O_W$ is admissible, so $O$-twisted $\C$-enriched categories with $W$ as the set of objects and functors that induce identity on objects carry a model structure
whose weak equivalences and fibrations are those $\C$-enriched functors $F\colon\D\to\mathcal E$ that induce weak equivalences, respectively, fibrations in~$\C$:
$$\Hom_\D(D,D')\to\Hom_\E(C,C'),$$
for all objects $D=F(D)$~and~$D'=F(D')$ in~$\Ob(\D)=\Ob(\mathcal E)=W$.

If $\C$ is in addition flat over the levels $\varphi_n$ ($n \ge 0$) of some weak equivalence $\varphi\colon O \r P$ of nonsymmetric operads in~$\V$, there is a Quillen equivalence of $O$- and $P$-twisted $\C$-enriched categories (both with $W$ as objects):
$$\varphi_* : \Cat^O_W(\C) \leftrightarrows \Cat^P_W(\C): \varphi^*.$$
For example, if $1_\V$ is cofibrant, then this condition is satisfied for any weak equivalence
$\Ai\to\Ass$, where $\Ai$ is a cofibrant replacement of $\Ass$.
It is satisfied for any weak equivalence if $\C$ is flat (\refde{summary.symmetricity}).
\xcoro

\pf
Admissibility follows from~\refth{O.Alg} and \refre{module.context} and rectification follows from~\refth{rect.colored.operad}.
If $1_\V$ is cofibrant, then $\C$ is flat over the levels of $\Ai \r \Ass$: $\Ass_n = 1_\V$ is cofibrant.
Moreover, $\Ai$ is a cofibrant operad, so that its levels are cofibrant by \refle{cofibrancy.preserved}.
Any monoidal model category is flat over a weak equivalence between cofibrant objects by Brown's lemma.
\xpf

These individual model structures on~$\Cat_W(\C)$ can be assembled into a single model structure on~$\Cat(\C)$.
The following result is due to Muro~\cite[Theorem~1.1]{Muro:Dwyer}, for $\C$ combinatorial.
Muro's work relaxes the assumptions of a similar result
of Berger and Moerdijk~\cite[Theorem~1.9]{BergerMoerdijk:Homotopy},
which in turn generalizes results of
Amrani ($\V = \Top$) \cite{Amrani:Model},
Bergner (for $\V = \sSet$) \cite[Theorem~1.1]{Bergner:Model},
Lurie (every object of~$\V$ is cofibrant) \cite[Proposition~A.3.2.4]{Lurie:HTT},
and
Tabuada ($\V = \Ch(\Mod_R)$ for some ring $R$ and $\V$ being symmetric spectra) \cite[Th\'eor\`eme~3.1]{Tabuada:Invariants}, \cite{Tabuada:Erratum}, and \cite[Theorem~5.10]{Tabuada:Homotopy}.

Given some property of objects or morphisms in~$\C$, we say that a $\C$-enriched category or a $\C$-enriched functor has this property
locally if it is true for the enriched objects of morphisms between each pair of objects.
Given a $\C$-enriched category, its derived~$\pi_0$ is an ordinary 1-category that is constructed by applying the derived internal hom from the monoidal unit of~$\C$
to each object of morphisms.

\prop (Muro)
\label{prop--Cat.C}
Suppose again that $\C$ is h-monoidal, and moreover combinatorial.
Then $\Cat(\C)$ carries the \emph{Dwyer--Kan model structure} whose weak equivalences are the \emph{Dwyer--Kan equivalences}
(i.e., local weak equivalences and their derived~$\pi_0$ is an essentially surjective functor or, equivalently, an equivalence of categories)
and whose acyclic fibrations are local acyclic fibrations that are surjective on objects.
\xprop

\rema
We expect that the preceding result can be extended to the case when $\C$ is not combinatorial,
but just admissibly generated and tractable using \refth{O.Alg}\refit{admissibly.generated.case}.
More generally, we expect that the Dwyer--Kan model structure on $\Cat^O(\C)$ exists for any operad $O$.
The reader is encouraged to upgrade Muro's result \ref{prop--Cat.C} to this generality, which will give a model structure on topological categories.
\xrema

\prop
Fix $\V$ and $\C$ as in \refco{Cat.C.fixed.objects} and a weak equivalence $\varphi:O\to P$ of nonsymmetric operads in~$\V$.
Assume that the Dwyer--Kan model structure on $\Cat^O(\C)$ and $\Cat^P(\C)$ exists, as in \refpr{Cat.C}.
If $\C$ is flat over a weak equivalence~$\varphi:O\to P$ (more precisely, flat over the levels $\varphi_n$ for all $n \ge 0$),
then we have a Quillen equivalence
$$\varphi_* : \Cat^O(\C) \leftrightarrows \Cat^P(\C): \varphi^*.$$
For example, this holds for all weak equivalences $\varphi$ if $\C$ is flat.
It also holds for the weak equivalence $\varphi\colon \Ai \r \Ass$ if the monoidal unit $1_\V$ is cofibrant.
\xprop

\pf
For some cofibrant object $X \in \Cat^O(\C)$ and a fibrant object $Y \in \Cat^P(\C)$, the (co)unit morphism of the adjunction for $X$ and $Y$ can be computed in the corresponding slices $\Cat^O_{\Obj(X)}(\C)$ and
$\Cat^P_{\Obj(Y)}(\C)$.
Moreover, the (co)fibrancy of $X$ and $Y$ is equivalent to the one in the corresponding slice category.
Now the Quillen equivalence immediately follows from the rectification of category structures with a fixed set of objects (\refco{Cat.C.fixed.objects}).
\xpf

An interesting question that arises in relation to these results is whether it is possible
to define a monoidal structure on the category of enriched categories in such a way that
the resulting model category is monoidal.
The naive choice (take the product of sets of objects and the tensor product of enriched morphisms)
already fails to satisfy the pushout product axiom in the case when $\C$~is the model category of small categories,
as shown by Lack.
The {\it Gray tensor product\/} does turn enriched categories in small categories (i.e., strict 2-categories)
into a monoidal model category, however, it is unclear how one should generalize it to enriched categories.
If such a monoidal product could be constructed, then one could iterate the construction of enriched categories
and consider higher enriched categories (i.e., enriched categories in enriched categories, etc.).
Such a construction could explain how the traditional definitions of bicategories, tricategories, and tetracategories
could be generalized in a systematic way to higher dimensions.
Furthermore, for certain choices of the operad~$O$ (e.g., the categorical $\Ai$-operad) one would expect
to get a model category that is Quillen equivalent to any of the usual model categories of $(\infty,n)$-categories.
(We cannot expect this for $O=\Ass$ because it is well known that tricategories cannot in general be strictified to strict 3-categories.)

\subsection{Applications to category theory}
\label{sect--category.theory}

In this section, we apply the results of \refsect{enriched.categories} to some concrete examples of (low-dimensional) category theory.

Consider the category of sets equipped with the model structure whose weak equivalences are bijections
and fibrations and cofibrations are arbitrary maps.
Equip this model category with the monoidal structure given by the cartesian product.
This model structure is tractable, proper,
its weak equivalences are stable under filtered colimits (it is pretty small in the sense of \csy{\refde{pretty.small}} for the maps $\emptyset\to\{0\}$, $\{0,1\}\to\{0\}$ generate the cofibrations, then use \csy{\refle{sequential}}),
symmetric h-monoidal and symmetroidal, and symmetric flat.
By \refpr{Cat.C}, the category $\Cat$ of categories admits a model structure whose weak equivalences
are equivalences of categories and fibrations are the so-called isofibrations,
i.e., functors $F\colon \C\to\D$ such that any isomorphism in~$\D$, $F(C) \cong D$ (for $C \in \C$, $D \in \D$) has a lift to an isomorphism in $\C$.
This is precisely the {\it canonical (folk) model structure on categories\/} (see, for example, Rezk~\cite{Rezk:Model}).
The canonical model structure is tractable, pretty small, cartesian (i.e., monoidal with
respect to the categorical product), simplicial, and all objects are fibrant and cofibrant (see Rezk~\cite{Rezk:Model} for details).
Furthermore, it is symmetric h-monoidal and symmetroidal because cofibrations are precisely those functors that are injective on objects,
and the latter property survives pushout products and coinvariants under~$\Sigma_n$, the argument being similar to the one for simplicial sets (see \csy{\refsect{simplicial.sets}}).
Finally, the canonical model structure is flat, which follows immediately from the definition of equivalences of categories, which are stable under products.
However, symmetric flatness fails: the $\Sigma_n$-equivariant functor from the groupoid~$\E\Sigma_n$ (objects are~$\Sigma_n$
and morphisms are $\Sigma_n\times\Sigma_n$) to the terminal groupoid is a weak equivalence,
yet its $\Sigma_n$-coinvariants is the map $\B\Sigma_n\to1$ ($\B\Sigma_n$ has one object whose endomorphisms are~$\Sigma_n$),
which is not an equivalence.

The results of \refsect{operad.admissible}--\refsect{rectification} yield model structures on various types of monoidal categories and a strong form of
\emph{Mac~Lane's coherence theorem}.

\prop
There is a model structure on strict monoidal categories, monoidal categories,
strict symmetric monoidal categories, and symmetric monoidal categories whose weak equivalences and fibrations are the ones of the underlying categories.

Every monoidal category is equivalent (via a strong monoidal functor) to a strict monoidal category.
This strict monoidal category is unique up to strict monoidal equivalence.
Similarly, every monoidal functor is equivalent (via a strong monoidal natural transformation) to a strict monoidal functor,
which is again unique up to a strict monoidal natural transformation.
\xprop

\pf
The above-mentioned categories are algebras (in $\Cat$) over the associative operad~$\Ass$, operad~$\Ai$, commutative operad~$\Comm$, and operad~$\Ei$, respectively.
Hence, the existence of the model structure follows from \refth{O.Alg}, whose assumptions have been verified above.

Furthermore, the nonsymmetric rectification theorem (\refth{rect.colored.operad}) tells us that the canonical morphism from~$\Ai$ to the associative operad induces
a Quillen equivalence between $\Ass$-algebras and $\Ai$-algebras.
\xpf

\exam
The morphism from~$\Ei$ to the commutative operad is not symmetric flat, as explained above, which tells us that symmetric monoidal
categories cannot always be strictified to strict symmetric monoidal categories.
This is well known
because symmetric monoidal categories can have a nontrivial k-invariant whereas strict symmetric monoidal categories always have a trivial k-invariant.
\xexam

Similarly, \emph{Mac~Lane's coherence theorem for bicategories} follows from the above, since
strict 2-categories are $\CatOp^\Ass$-algebras and bicategories are $\CatOp^\Ai$-algebras in $\Cat$, respectively.

\prop
There is a Quillen equivalence between the model categories of
strict 2-categories and bicategories.
\xprop

We conjecture that other strictification results of category theory,
such as strictification of tricategories to Gray categories (Gordon, Power, and Street),
partial strictification of symmetric monoidal bicategories, etc., can also be shown using the methods of this paper.
However, considerations of volume prevent us from developing this topic further.
Simpson's conjecture might also be amenable to the techniques explained above.

\subsection{The colored operad of colored operads}
\label{sect--operad.operads}
Given a set~$W$, there is a (symmetric) colored operad~$\OpOp_W$ whose category of algebras
is equivalent to the category of (symmetric) $W$-colored operads in $\C$.
It is due to Berger and Moerdijk \cite[\S1.5.6 and \S1.5.7]{BergerMoerdijk:Resolution}.
See also \cite[\S3]{GutierrezVogt:Model} for a detailed description of the multicolored case.

This operad is first constructed for $\C = \Sets$ as follows: the set of colors of~$\ssOpOp_W$ is the set of objects of~$\ssSeq_{W,W}$, which we call {\it valencies}.
Recall from \refsect{colored.collections} that the objects of $\ssSeq_{W,W}$ are pairs $c = (s,w)$, where $s\colon I \r W$ is a map from a finite set~$I$ and $w \in W$.
The operations
$$\ssOpOp_W (a_1,\ldots,a_k;b)$$ from a given sequence of valencies~$(a_1,\ldots,a_k)$ to a valency~$b$
are given by isomorphism classes of triples $(T, \sigma, \tau)$ consisting of a $W$-colored (symmetric) tree~$T$ equipped with a bijection~$\sigma$ from~$\{1,\ldots,k\}$ to the set of internal vertices of~$T$
such that the valency of~$\sigma(i)$ equals~$a_i$
and a color-preserving bijection~$\tau$ from $\{1,\ldots,m\}$, where $m$ is the arity of~$b$, to the input edges of~$T$.
Isomorphisms of such triples are isomorphisms of colored trees that are compatible with $\sigma$ and $\tau$.
In the symmetric case, the symmetric group~$\Sigma_k$ acts on such classes by precomposition with~$\sigma$.
The operadic unit sends each valency~$c$ to the corresponding corolla, interpreted as an operation from~$c$ to~$c$.
The operadic composition is given by inserting trees into each other (see \cite[\S1.5.6]{BergerMoerdijk:Resolution} in the uncolored case).
One checks that this gives a (symmetric) operad, denoted by $\ssOpOp_W$, in $\Sets$.

The functor
$$\Sets \r \C, \qquad X \mapsto \coprod_{x \in X} 1_\C$$
is symmetric monoidal and therefore extends to a functor
$$\ssOper_{\ssSeq_W}(\Sets) \r \ssOper_{\ssSeq_W}(\C).$$
The image of $\ssOpOp_W$ under this functor is again denoted by $\ssOpOp_W$.

The following admissibility statement unifies a few earlier results: the semimodel structure for symmetric operads established by Spitzweck \cite[Theorem~3.2]{Spitzweck:Operads}, the model structure for nonsymmetric operads by Muro \cite[Theorem~1.1]{Muro:Homotopy}, and the model structure on uncolored operads in orthogonal spectra with the positive stable model structure by Kro \cite[Theorem~1.1]{Kro:Model}.

\coro \label{coro--operads}
Let $\C$ be (symmetric) h-monoidal.
Then the operad $\ssOpOp_W$ of (symmetric) $W$-colored operads is admissible, that is to say, the category~$\ssOper_W(\C)$ of (symmetric) $W$-colored operads in~$\C$
has a model structure that is transferred along the adjunction
$$\Free: \C^{\ssSeq_{W,W}}\rightleftarrows\Alg_{\OpOp_W}(\C)=\ssOper_W(\C): U.$$
If $1_\C$ is cofibrant, then $\ssOpOp_W$ is strongly admissible, i.e., the forgetful functor~$U$ preserves cofibrations with cofibrant domain.
\xcoro

\pf
The admissibility follows from \refth{O.Alg}.
The strong admissibility follows from \refpr{cofibrant.strongly.admissible} since $\ssOpOp$ is levelwise projectively cofibrant.
\xpf

Enriched operads are generalized to enriched $\Ai$-operads in the same fashion as enriched categories are generalized to enriched $\Ai$-categories.
Fix a (symmetric) operad~$O$.
In practice, $O$ is an $\Ai$-operad, i.e., we have a weak equivalence of operads $O\to\Ass$, where $\Ass$ denotes the associative operad.
We define the colored (symmetric) operad~$\OpOp^O_W$ of \emph{$O$-twisted $W$-colored (symmetric) operads} by the same construction as above,
starting from a colored operads~$P$ in sets,
except that we pass to a $\C$-valued operad in a modified fashion:
instead of tensoring operations in degree~$k$ with $1_\C$, we tensor them with $O_k$.
The intuitive idea behind this is that the composition of operadic operations is no longer strictly associative,
but is rather governed by the operad~$O$.
An $O$-twisted $W$-colored operad is an $O$-algebra in the monoidal category of $W$-colored (symmetric) sequences equipped with the substitution product,
the latter being a left $\C$-module in the obvious way.
Then \refco{operads} has an immediate generalization for the operad $\ssOpOp^O$.
For the strong admissibility, the requirement on $1_\C$ is replaced by the condition that the levels $O_k$ be cofibrant as objects in $\C$.
Moreover, \refth{rect.colored.operad} admits the following corollary.

\coro \label{coro--operads.2}
If $\C$ is flat over a weak equivalence~$O\to P$ of operads, then we have a Quillen equivalence $\sOper^O_W(\C)\leftrightarrows\sOper^P_W(\C)$ of $O$-twisted and $P$-twisted
(symmetric) $W$-colored operads in~$\C$.
For example, if $1_\C$ is cofibrant, then $\Ai$-twisted colored symmetric operads can be rectified to ordinary colored symmetric operads.
\xcoro

\pf
This follows from \refth{rect.colored.operad} once we show the symmetric flatness of~$\C$ with respect to $\sOper^O_W\to\sOper^P_W$.
Every component of $\sOper^O_W$ is a coproduct of the corresponding components of~$O$,
and the relevant symmetric group acts freely on the components.
Thus, the symmetric flatness follows from the flatness of~$\C$ over $O\to P$.
\xpf

\rema
In fact, if $\C$ is a $\V$-enriched model category that is symmetric h-monoidal with respect to~$\V$ only (and not necessarily with respect to itself),
then the colored operad of colored operads can be defined with values in~$\V$ and its algebras in~$\C$ will still be $W$-colored operads in~$\C$,
so the above corollary holds in this more general setting.
Guti\'errez and Vogt used such a setup (with a different set of conditions on~$\V$)
to construct a model structure on $W$-colored operads in symmetric spectra (see \cite[Corollary~4.1]{GutierrezVogt:Model}).
\xrema

Starting from this point, further work is required to assemble the model structures on~$\sOper_W(\C)$
into one on the category $\ssOper(\C)$ of (symmetric) operads with an arbitrary set of colors.
This has been done for $\C=\sSet$ by Cisinski and Moerdijk \cite[Theorem~1.14]{CisinskiMoerdijk:Dendroidal}
and independently by Robertson~\cite[Theorem~6]{Robertson:Homotopy}
and was extended by Caviglia~\cite{Caviglia:Model} to more general model categories using similar arguments.
We expect that the assumptions can be further relaxed to the ones stated in the above corollary.

\subsection{Diagrams}

In this section, we construct a model structure on the category of enriched diagrams of some fixed shape and prove a rectification result.
In particular, we recover the classical result of Vogt and its generalization by Cordier and Porter
on homotopy coherent diagrams.

\prop
Assume that $\C$ is, in addition to the standing assumptions in this section, h-monoidal.
For any $\V$-enriched, small category $\D$,
the category of $\V$-enriched functors $\D\to\C$ admits a transferred model structure.
Its weak equivalences and fibrations are those natural transformations of $\V$-enriched functors $F \r G$
such that for all objects $X \in \D$,
$$F(X) \r G(X)$$
is a weak equivalence, respectively, a fibration.
Furthermore, if $\V$ has a model structure and $\C$ is flat over~$\V$,
then a componentwise weak equivalence of diagrams $\D\to\D'$ whose object map is the identity
induces a Quillen equivalence of the two model categories of diagrams.
\xprop

\rema
A more general version of the rectification result allows for a Dwyer--Kan equivalence $\D\to\D'$.
\xrema

\pf
Following Berger and Moerdijk~\cite[\S1.5.5]{BergerMoerdijk:Resolution},
we consider the nonsymmetric colored operad $\DiagOp_\D$
that encodes diagrams in~$\C$ indexed by a fixed $\V$-enriched category~$\D$,
i.e., $\V$-enriched functors $\D\to\C$.
The operad $\DiagOp_\D$ is colored by the set of objects of~$\D$.
Its operations are defined as
$$\DiagOp_\D(X_1, \ldots, X_n, Y)=
\cases{
\emptyset, & $n \ne 1$; \cr
\Map_\D(X, Y), & $n = 1$.}$$
Here, $\Map_\D$ denotes the enriched hom object.
The operadic composition and unit are induced by the composition and unit of~$\D$.
(The construction just described embeds enriched categories into nonsymmetric colored operads.)

A $\DiagOp_\D$-algebra in~$\C$ consists of a collection of objects~$D_X$ in~$\C$, for all $X \in \D$
together with morphisms $\Mor(X,Y)\otimes D_X\to D_Y$ that satisfy the obvious associativity and unitality conditions.
This is precisely the data of a $\V$-enriched functor $\D\to\C$.

\refth{O.Alg} now implies that the category of $\D$-diagrams admits a transferred model structure.
At this point, we remark that
\refth{O.strongly.admissible} likewise implies that cofibrations with cofibrant source are preserved by the forgetful functor
if taking the pushout product with $1_C\to\Map_\D(X,X)$ and $\emptyset\to\Map_\D(X,Y)$ preserves (acyclic) cofibrations,
which is true, for example, if individual hom objects are cofibrant and the unit maps are cofibrations.

\refth{rect.colored.operad} implies the desired rectification statement if $\C$ is flat.
\xpf

\subsection{Monoidal diagrams}

Extending the results of the previous section, there is also a (symmetric) colored operad that encodes lax (symmetric) monoidal diagrams,
i.e., lax (symmetric) monoidal $\V$-enriched functors $\D\to\C$, where $\C$ is now an algebra over the monoidal category~$\V$
and $\D$ is a monoidal $\V$-enriched category.
We therefore obtain a model structure on lax (symmetric) monoidal functors.

\prop
Assume that $\C$ is (symmetric) h-monoidal.
For any $\V$-enriched symmetric monoidal small category $\D$,
the category of lax (symmetric) monoidal $\V$-en\-riched functors $\D\to\C$ admits a transferred model structure.
Furthermore, if $\C$ is (symmetric) flat over~$\V$, then a weak equivalence $\D\to\D'$ induces a Quillen equivalence of the induced model categories.
\xprop

\pf
We consider the (symmetric) operad whose operations from a multisource $(s_1,\ldots,s_k)$ to a target~$t$
are given by the enriched morphism object from~$s_1\otimes\cdots\otimes s_k$ to~$t$.
The operadic composition and unit are induced by the monoidal category structure of~$\D$.

An algebra in~$\C$ over this operad consists of a collection of objects~$D_X$ in~$\C$, for any $X \in \D$,
together with morphisms $\Mor(X_1\otimes\cdots\otimes X_k,Y)\otimes D_{X_1}\otimes\cdots\otimes D_{X_k}\to D_Y$ that satisfy the corresponding associativity and unitality conditions.
This is precisely the data of a (symmetric) lax monoidal $\V$-enriched functor $\D\to\C$.

As before, Theorems \ref{theo--O.Alg} and \ref{theo--rect.colored.operad}
now imply the admissibility and rectification criteria as stated.
\xpf


\subsection{Prefactorization algebras}

As an application of the previous section, we construct a model structure on prefactorization algebras.
See Costello and Gwilliam's book~\cite[\S7.3]{CostelloGwilliam:Factorization} for the relevant background.
A \emph{prefactorization algebra} on a $\V$-enriched monoidal site~$(S, \sqcup, \emptyset)$
(it is useful to think of the monoidal structure as the disjoint union)
is a symmetric lax monoidal $\V$-enriched functor from~$S$ to~$\C$, where $\C$ is $\V$-enriched.
A typical example of~$S$ is the category of smooth manifolds and their embeddings equipped with the Weiss topology,
where morphism objects are either discrete or have the natural space structure.
The previous section now immediately implies the following statement.

\prop
If $\C$ is symmetric h-monoidal, $\V$-enriched, and $S$ is a $\V$-enriched site, then the category of prefactorization algebras over~$S$ with values in~$\C$ admits a transferred model structure.
Furthermore, if $\C$ is symmetric flat, then a functor of sites $S\to S'$ that induces the identity morphism on objects
and is a componentwise weak equivalence on morphism
gives a Quillen equivalence of the corresponding model categories.
\xprop

This raises the question whether the above model structure can be upgraded to factorization algebras.

\raggedright\rightskip0em plus \maxdimen

\bibliographystyle{dp}
\def\ZM#1{\href{https://zbmath.org/?q=an:#1}{Zbl #1}}
\def\MRx#1 #2\relax{\href{https://mathscinet.ams.org/mathscinet-getitem?mr=#1}{MR#1}}
\def\MR#1{\MRx#1 \relax}
\bibliography{bib}

\end{document}